\documentclass[12pt]{amsart}
\usepackage{amsmath,amscd,amssymb,amsfonts,graphics,mathrsfs}
\usepackage{hyperref}
\newcommand{\hl}{\hyperlink}
\newcommand{\htt}{\hypertarget}
\newcommand{\h}{\hbox}
\newcommand{\q}{\quad}
\newcommand{\nin}{\noindent}
\newcommand{\bs}{\par\bigskip}
\newcommand{\ms}{\par\medskip}
\newcommand{\sk}{\par\smallskip}
\newcommand{\bsn}{\par\bigskip\noindent}
\newcommand{\msn}{\par\medskip\noindent}
\newcommand{\skn}{\par\smallskip\noindent}
\newcommand{\ges}{\geqslant}
\newcommand{\les}{\leqslant}
\newcommand{\1}{\hskip1pt}

\newcommand{\mcup}{\hbox{$\bigcup$}}

\newcommand{\msum}{\hbox{$\sum$}}
\newcommand{\mprod}{\hbox{$\prod$}}

\newcommand{\D}{{\mathscr D}}
\newcommand{\G}{{\mathscr G}}

\newcommand{\bt}{\widetilde{b}}
\newcommand{\ct}{\widetilde{c}}
\newcommand{\Ft}{\widetilde{F}}

\newcommand{\Rt}{\widetilde{R}}

\newcommand{\alt}{\widetilde{\alpha}}
\newcommand{\nut}{\widetilde{\nu}}

\newcommand{\R}{{\mathbb R}}

\newcommand{\Q}{{\mathbb Q}}
\newcommand{\C}{{\mathbb C}}
\newcommand{\N}{{\mathbb N}}

\newcommand{\Z}{{\mathbb Z}}

\newcommand{\Hppf}{\mathscr{H}''_{\!f}}
\newcommand{\Htppf}{\widetilde{\mathscr{H}}''_{\!f}}
\newcommand{\delf}{\delta_{\!f}}

\newcommand{\Gf}{\G_f}
\newcommand{\dti}{\dd_t^{-1}}
\newcommand{\Cdti}{\C\{\!\{\dd_t^{-1}\}\!\}}
\newcommand{\ob}{{\bf 1}}
\newcommand{\Gr}{{\rm Gr}}

\newcommand{\al}{\alpha}
\newcommand{\be}{\beta}
\newcommand{\ep}{\varepsilon}
\newcommand{\ga}{\gamma}
\newcommand{\la}{\lambda}
\newcommand{\om}{\omega}
\newcommand{\Om}{\Omega}

\newcommand{\dd}{\partial}
\newcommand{\ddd}{{\rm d}}
\newcommand{\eq}{\,{=}\,}
\newcommand{\gess}{\,{\ges}\,}
\newcommand{\less}{\,{\les}\,}
\newcommand{\sgt}{\,{>}\,}
\newcommand{\slt}{\,{<}\,}
\newcommand{\nes}{\,{\ne}\,}
\newcommand{\notins}{\,{\notin}\,}
\newcommand{\mi}{\1{-}\1}

\newcommand{\pl}{\1{+}\1}
\newcommand{\scd }{\1{\cdot}\1}
\newcommand{\bl}{\bigl}
\newcommand{\br}{\bigr}
\newcommand{\sst}{\,{\subset}\,}
\newcommand{\stm}{\,{\setminus}\,}
\newcommand{\ins}{\,{\in}\,}
\newcommand{\tos}{\,{\to}\,}
\newcommand{\defs}{\,{:=}\,}
\newcommand{\ssc}{\,\raise.15ex\hbox{${\scriptstyle\circ}$}\,}
\newcommand{\ssb}{\raise.15ex\h{${\scriptscriptstyle\bullet}$}}
\newcommand{\into}{\hookrightarrow}

\newcommand{\simto}{\,\,\rlap{\hskip1.3mm\raise1.4mm\hbox{$\sim$}}\hbox{$\longrightarrow$}\,\,}

\begin{document}
\h{}\bs
\centerline{\large Bernstein-Sato Polynomials of Semi-weighted-homogeneous}{
\sk
\centerline{\large Polynomials of nearly Brieskorn-Pham Type}
\bs
\centerline{Morihiko Saito}
\bsn\ms
\vbox{\nin\narrower\smaller
{\bf Abstract.} Let $f$ be a semi-weighted-homogeneous polynomial having an isolated singularity at 0. Let $\alpha_{f,k}$ be the spectral numbers of $f$ at 0. By Malgrange and Varchenko there are non-negative integers $r_k$ such that the $\alpha_{f,k}{-}\1r_k$ are the roots up to sign of the local Bernstein-Sato polynomial $b_f(s)$ divided by $s{+}1$. However, it is quite difficult to determine these shifts $r_k$ explicitly on the parameter space of $\mu$-constant deformation of a weighted homogeneous polynomial. Assuming the latter is nearly Brieskorn-Pham type, we can obtain a very simple algorithm to determine these shifts, which can be realized by using Singular (or even C) without employing Gr\"obner bases. This implies a refinement of classical work of M. Kato and P. Cassou-Nogu\`es in two variable cases, showing that the stratification of the parameter space can be controlled by using the (partial) additive semigroup structure of the weights of parameters. As a corollary we get for instance a sufficient condition for all the shiftable roots of $b_f(s)$ to be shifted. We can also produce examples where the minimal root of $b_f(s)$ is quite distant from the others as well as examples of semi-homogeneous polynomials with roots of $b_f(s)$ nonconsecutive.}
\ms\bs
\centerline{\bf Introduction}
\bsn
Let $f\ins\C\{x\}$ be a convergent power series of $n$ variables having an {\it isolated\1} singularity at 0, where $f(0)\eq0$ and $n\gess 2$. Set $(X,0)\defs(\C^n,0)$. The local Bernstein-Sato polynomial $b_f(s)$ (see \cite{Be}, \cite{SaSh}, \cite{Bj}, \cite{Sat}, \cite{Ka}) is called the {\it BS polynomial\1} for short in this paper. (This is a factor of the {\it global\1} BS polynomial in the $f$ polynomial case using finite determinacy of holomorphic functions with isolated singularities, since $f$ may have singularities outside the origin.) Let $\Rt_f\sst\R_{>0}$ be the set of roots up to sign of the {\it reduced\1} BS polynomial $\bt_f(s)\defs b_f(s)/(s{+}1)$, see \cite{Ka} for negativity of the roots of $b_f(s)$. Let $\al_{f,k}$ ($k\ins[1,\mu_f]$) be the {\it spectral numbers\1} of $f$ counted with multiplicities, where $\mu_f$ is the Milnor number, see \cite{St2} (and also \cite{DS}). We assume that the $\al_{f,k}$ are weakly increasing. As a consequence of theorems of Malgrange \cite{Ma} and Varchenko \cite{Va1}, it is well known to specialists that there are {\it non-negative integers\1} $r_k$ ($k\ins[1,\mu_f]$) such that
\htt{1}{}
$$\Rt_f\eq\{\al_{f,k}{-}\1r_k\}_k\q\h{with}\q\min\{\al_{f,k}\}_k\eq\min\{\al_{f,k}{-}\1r_k\}_k,
\leqno(1)$$
forgetting the multiplicities. We say that $\al_{f,k}{-}\1r_k$ is a {\it shifted\1} root up to sign of the reduced BS polynomial $\bt_f(s)$ if $r_k\gess 1$. It is, however, quite nontrivial to determine these shifts $r_k$ on the base space of the miniversal $\mu$-constant deformation of a weighted homogeneous polynomial. This is a very interesting problem, see for instance \cite{Kat1}, \cite{Kat2}, \cite{Ca}, \cite{NT}. (Recall that the spectral numbers are {\it invariant\1} by $\mu$-constant deformations, and coincide with the roots of $\bt_f(s)$ up to sign in the {\it weighted homogeneous\1} case.)
\sk
The above assertion (\hl{1}{1}) is easily recognized if we define the {\it saturated Hodge filtration\1} $\Ft$ on the {\it $\la$-eigenspaces of the vanishing cohomology\1} $H^{n-1}(F_{\!f},\C)_{\la}\defs{\rm Ker}(T_s{-}\la)$ (with $F_{\!f}$ the Milnor fiber and $T\eq T_sT_u$ the Jordan decomposition of the monodromy) replacing the {\it Brieskorn lattice\1}
\vskip-7mm
$$\Hppf\defs\Om^n_{X,0}/\ddd f{\wedge}\ddd\Om^{n-2}_{X,0},$$
(see \cite{Br}) by its {\it saturation\1} $\Htppf\defs\msum_{i\ges0}\,(\dd_tt)^i\Hppf$ (see \cite{Ma}) in the formula for the Hodge filtration $F$ in \cite[(4.5)]{SS} (or \cite[(2.6.3)]{bl}, see (\hl{1.1.8}{1.1.8}) below) and considering the meaning of the bigraded pieces $\Gr_{\Ft}^{p+r}\Gr_F^pH^{n-1}(F_{\!f},\C)_{\la}$, since $r$ corresponds to the {\it shift.} Here $H^{n-1}(F_{\!f},\C)_{\la}$ is identified with $\Gr_V^{\al}\Gf$ for $\la\eq e^{-2\pi i\al}$ with $\Gf\defs \Hppf[\dd_t]$ the localization of $\Hppf$ by the action of $\dti$, which is called the {\it Gauss-Manin system\1} (see for instance \cite{bl}), and $V$ denotes the filtration of Kashiwara and Malgrange on the regular holonomic $\D_{\C,0}$-module $\Gf$ {\it indexed by\1} $\Q$. (This filtration was originally indexed by $\Z$, see also \cite{sup}, \cite[\S3.4]{fil} about the reason for which $V$ must be indexed by $\Q$.)
\sk
We assume in the introduction that the monodromy $T$ is {\it semisimple,} that is, $N\defs\log T_u$ vanishes. Let $\om_k$ ($k\ins[1,\mu_f]$) be free generators of the Brieskorn lattice $\Hppf$ over $\Cdti$ (see for instance \cite{bl} for $\Cdti$). We {\it assume\1} that the $\om_k$ give a $\C$-basis of the quotient space
\vskip-7mm
$$\Omega_f^n\defs\Om^n_{X,0}/\ddd f{\wedge}\Om^{n-1}_{X,0}\eq\Hppf/\dti\Hppf,$$
in a {\it compatible way\1} with the $V$-filtration (that is, inducing a $\C$-basis of $\Gr_V^{\ssb}\Omega_f^n$); for instance, the $\om_k$ are associated with an {\it opposite filtration\1} in the sense of \cite{bl}. Using the {\it semi-simplicity\1} of $T$, we can set
$$\Gf^{(\al)}\defs{\rm Ker}(\dd_tt\mi\al)\sst\Gf\q(\al\in\Q),$$
so that we have the {\it asymptotic expansions\1}
\htt{2}{}
$$\om_k\eq\msum_{\al\ges\al_{f,k}}\,\om_k^{(\al)}\q(k\ins[1,\mu_f]),
\leqno(2)$$
with $\om_k^{(\al)}\ins\Gf^{(\al)}$ and $v_k\defs\om_k^{(\al_{f,k})}\nes 0$ (renumbering the $\om_k$ if necessary). The $v_k$ form a free basis of $\Gf$ over $\Cdti[\dd_t]$, and we have the {\it power series expansions}
\htt{3}{}
$$\om_k\eq\msum_{l=1}^{\mu_f}\,g_{k,l}\1v_{l}\q\h{with}\q g_{k,l}\ins\Cdti[\dd_t].
\leqno(3)$$
Assume furthermore the $\om_k$ give a $\C$-basis of $\Gr_{\Ft}^{\ssb}\Gr_F^{\ssb}H^{n-1}(F_{\!f},\C)_{\la}$ (inducing a {\it bisplitting\1} of $F,\Ft$), where the above identification of $H^{n-1}(F_{\!f},\C)_{\la}$ with $\Gr_V^{\al}\Gf$ is used. This condition is not always satisfied in the semi-weighted-homogeneous case if we consider only {\it monomial bases.} It is, however, trivially satisfied if the spectral numbers have {\it multiplicity\1} 1, that is,
\htt{M1}{}
$$\al_{f,k}\nes\al_{f,l}\q(k\nes l).
\leqno{\rm(M1)}$$
The last condition holds if every eigenvalue of the monodromy has {\it multiplicity\1} 1, that is,
\htt{M1'}{}
$$\dim H^{n-1}(F_{\!f},\C)_{\la}\less 1\q(\forall\,\la),
\leqno{\rm(M1)}{}'$$
for instance, if $f$ is a Brieskorn-Pham type polynomial with {\it mutually prime} exponents. The above argument is not appropriate for explicit calculations of $\Htppf$ unless condition~(\hl{M1}{M1}) is satisfied at least partially, since it is not easy to find free generators $\om_k$ satisfying the condition related to $\Gr_{\Ft}^{\ssb}\Gr_F^{\ssb}$ written after (\hl{3}{3}) {\it before determining the saturation\1} $\Htppf$ in general.
\sk
We can easily verify the following.
\par\htt{P1}{}\msn
{\bf Proposition~1.} {\it Under the above assumptions, the saturation $\Htppf$ is generated over $\Cdti$ by the $\om_k^{(\al)}$ for $k\ins[1,\mu_f]$, $\al\in[\al_{f,k},n\mi\alt_f]$, where $\alt_f\eq\al_{f,1}$ $($the minimal spectral number$)$.}
\ms
Combining this with the expansions in (\hl{3}{3}), we get the following.
\par\htt{C1}{}\msn
{\bf Corollary~1.} {\it With the above notation and assumptions, let $m_{k,l}$ be the pole order of $g_{k,l}$ as a power series of $\tau\defs\dti$ with a pole of finite order. Set $r_l\defs\max\{m_{k,l}\}_{k\in[1,\mu_f]}$ $(l\ins[1,\mu_f])$. Then the $\al_{f,k}\mi r_k$ $(k\ins[1,\mu_f])$ are the roots of $\bt_f(s)$ up to sign forgetting the multiplicities.}
\ms
Assume that $f$ is a {\it semi-weighted-homogeneous\1} polynomial $\msum_{\be\ges 1}\,f_{\be}$, where the $f_{\be}$ are weighted homogeneous polynomials of weighted degree $\be$, which vanish except for a finite number of $\be$. (Recall that convergent power series with isolated singularities have {\it finite determinacy,} see for instance \cite{GLS}.) We assume that the lowest weighted degree part $f_1$ is a polynomial of {\it Brieskorn-Pham type\1} (abbreviated as BP type) $\msum_{i=1}^n\,x_i^{e_i}$ with exponents $e_i\gess 3$; in particular, $f$ has an isolated singularity at 0, and the weight $w_i$ of $x_i$ is $1/e_i$.
\sk
By \cite{St1}, \cite{Va2}, the {\it spectral numbers\1} $\al_{f,k}$ ($k\ins[1,\mu_f]$) are given by
\htt{4}{}
$$\aligned&\Sigma_f=\bl\{\al_{\bf w}(\nu)\mid\nu\eq(\nu_i)\ins E_f\br\}\q\q\h{with}\\&\al_{\bf w}(\nu)\defs\msum_{i=1}^n\,\nu_i/e_i,\,\,\,\,E_f\defs \mprod_{i=1}^n\1\Z\,{\cap}\,[1,e_i{-}1].\endaligned
\leqno(4)$$
Set $\om^{\nu}\defs[x^{\nu-{\mathbf 1}}\ddd x]\ins\Hppf$ for $\nu\in\Z_{>0}^n$. Here $\ddd x\defs \ddd x_1{\wedge}\cdots{\wedge}\ddd x_n$ and ${\bf 1}\defs (1,\dots,1)$ (with ${\bf w}\eq (w_1,\dots,w_n)$). It is easy to see that the $\om^{\nu}$ for $\nu\ins E_f$ form a free basis of $\Hppf$ and we have the equalities $\al_{\bf w}(\nu)\eq\al_{\bf w}(\om^{\nu})$ ($\nu\ins E_f$) using an argument similar to \cite{exp}, where
$$\al_{\bf w}(\om)\defs\al_V(\om)\eq\max\{\al\ins\Q\mid\om\ins V^{\al}\Gf\}\q(\om\ins\Gf).$$
\sk
We now assume that condition~(\hl{M1}{M1}) holds, for instance, the $e_i$ are {\it mutually prime.} We may restrict the condition to $\al_{f,k}\notins\Z$ in the case $\alt_f\gess\tfrac{n}{2}{-}1$. Let $h_j$ ($j\ins J$) be the monomials such that $h_j$ is not contained in the Jacobian ideal $(\dd f_1)$ and moreover $\al_{\bf w}(h_j)\sgt1$. The first condition is equivalent to that $\nu{+}\ob\ins E_f$, and the second may be replaced with to the non-strict inequality $\al_{\bf w}(h_j)\gess1$ by condition~(\hl{M1}{M1}) assuming the first. We have the equality $\al_{\bf w}(h_j)\eq\al_{\bf w}(\nu)$ if $h_j\eq x^{\nu}$. So $\al_{\bf w}(h_j)\eq\al_{\bf w}([h_j\ddd x])\mi\alt_f$ with $\alt_f\eq\al_{f,1}\eq\al_{\bf w}([\ddd x])\eq\al_{\bf w}(\ob)$.
\sk
For the calculation of the BS polynomial, we may assume that $f$ is written as
\htt{5}{}
$$f=f_1\pl\msum_{j\in J}\,u_jh_j\q(u_j\ins\C),
\leqno(5)$$
using \cite{Va3} together with Remark~\hl{R1.3}{1.3} below. (The latter is related to a problem on the difference between {\it small\1} and {\it global\1} deformations.) Here $f_1$ is assumed to be BP type. In the case $n\eq2$, however, we can also allow that $f_1$ is {\it nearly Brieskorn-Pham type\1} (abbreviated as nearly BP type), that is, a linear combination of $n$~monomials having an isolated singularity at 0, which is called a non-degenerate invertible polynomial in mirror symmetry, see for instance \cite{Kre}, \cite{EG}. (Here $J$ must be modified slightly, see Remark~\hl{R1.1d}{1.1d} below.) Condition (\hl{M1}{M1}) implies that $|J|\eq {\rm mod}_{f_1}$, the {\it modality\1} of $f_1$. (The latter coincides with the inner modality by \cite{Va3}, see also \cite{NT}.) The $u_j$ are identified with the coordinates of the parameter space of the miniversal $\mu$-constant deformation of $f_1$, and have weights $\ga_j\defs\al_{\bf w}(h_j)\mi1$ for $j\ins J$. (Here the minus sign may be used in some papers.) We assume that the $\al_{f,k}$ and $\ga_j$ are {\it increasing\1} by identifying $J$ with $\{1,\dots,{\rm mod}_{f_1}\}$ so that
\htt{6}{}
$$\ga_j\eq\al_{f,k}\mi\alt_f\mi 1\,\,\,\,\h{if}\,\,\,\,k\eq j\pl\delf\,\,\,(j\ins J)\,\,\,\,\h{with}\,\,\,\,\delf\defs\mu_f\mi{\rm mod}_{f_1}.
\leqno(6)$$
Note that $\al_{f,k}\eq\al_{f_1,k}$ and $\alt_f\eq\alt_{f_1}$. We can calculate the Gauss-Manin connection quite explicitly, and get the following.
\par\htt{T1}{}\msn
{\bf Theorem~1.} {\it Let $f$ be a semi-weighted-homogeneous polynomial whose lowest weighted degree part $f_1$ is nearly BP type $($which is BP if $n\gess 3)$ as in {\rm(\hl{5}{5})}. Then there is an efficient algorithm $($without employing Gr\"obner bases$)$ to calculate the theoretically lowest coefficients of the $g_{k,l}\ins\Cdti[\dd_t]$ in \rm(\hl{3}{3})}}.
\ms
If $\alt_f\gess\tfrac{n}{2}{-}1$ (for instance, if $n\eq2$), we have either $r_{\!j}\eq1$ or 0 ($\forall\,j\ins J)$. It is then enough to calculate the {\it coefficients\1} $g_{k,l}^{(1)}$ of $\dd_t^1$ in the $g_{k,l}\ins\Cdti[\dd_t]$.
\par\htt{T2}{}\msn
{\bf Theorem~2.} {\it Assume $\alt_f\gess\tfrac{n}{2}{-}1$. With the notation and assumption of Theorem~{\rm\hl{T1}{1}}, the coefficients $g_{k,l}^{(1)}$ are weighted homogenous polynomials in the parameters $u_j$ $($with weight $\ga_j)$ of the miniversal $\mu$-constant deformation, and their weighted degrees are given by
\htt{7}{}
$$\deg_{\bf w}g_{k,l}^{(1)}\eq\deg_{\bf w}\om_{l}\mi\deg_{\bf w}\om_k\mi1\eq\al_{f,l}\mi\al_{f,k}\mi1.
\leqno(7)$$
Moreover the coefficient in $g_{k,l}^{(1)}$ of every monomial of the above weighted degree $($which is determined combinatorially$\1)$ is nonzero. Its sign depends only on whether the usual degree is even or odd $($that is, we get the positivity by replacing the $u_j$ with the $-u_j)$. The coefficient of $u_j$ in $g_{k,l}^{(1)}$ is $-1$ if $\ga_j$ coincides with the weighted degree of $g_{k,l}^{(1)}$.}
\ms
The last three assertions do not seem to follow easily from a general theory of Gauss-Manin connections. Note that $\deg_{\bf w}g_{1,l}^{(1)}\eq\ga_j\,({=}\,\al_{\bf w}(h_j)\mi1)$ in the case $\om_k\eq[\ddd x]$ (that is, $k\eq1$) and $\om_{l}\eq[h_j\ddd x]$ (that is, $l\eq j\pl\delf$). By (\hl{2}{2}) the weighted degrees $\deg_{\bf w}\om_k$ ($k\ins[1,\mu_f]$) coincide with the spectral numbers $\al_{f,k}$, and these are increasing.
\ms
It is not difficult to realize the above algorithm using Singular \cite{Sing} (or even C), see \hl{A1}{A.1}-\hl{A2}{2} in Appendix. The computation of necessary terms takes only a few seconds if the (usual) polynomial degree of $f_1$ is at most 8 in the two variable case. We have, however, a problem of {\it integer overflow\1} (even for $f_1\eq x^9{+}\1y^7$), and this may be eased {\it to some extent\1} by employing C which can treat 64 bit integers. This is rather a {\it structural problem,} since iterations of $\dd_t$ involve successive derivations of monomials of high degrees producing huge integers easily.
\sk
It is well known that there is a stratification of the parameter space $V$ of the miniversal $\mu$-constant deformation of $f_1$ such that the BS polynomial is constant on each stratum. This stratification, called the {\it BS stratification\1} in this paper, can be describes as follows: We identify $E_f$ with $\Z\,{\cap}\,[1,\mu_k]$ in such a way that the spectral numbers $\al_{f,k}$ ($k\ins[1,\mu_k]$) are increasing; in particular $\ob\ins\N^n$ corresponds to $1\ins[1,\mu_k]$. For $\nu,\nu^{\1\prime}\ins\N^n$, we say that $\nu^{\1\prime}$ (or $x^{\nu'}$) is {\it over\1} $\nu$ (or $x^{\nu}$), and note $\nu^{\1\prime}\succ\nu$, if $\nu^{\1\prime}_i\gess\nu_i$ $(\forall\,i)$, and similarly for {\it under\1} and $\prec$. For $j\ins J$, there is a unique $\nu^{(j)}\ins E_f$ with $x^{\nu^{(j)}}\eq h_j$ (assuming condition~{\rm(\hl{M1}{M1})} for $\al_{f,k}\notins\Z$ with $\alt_f\gess\tfrac{n}{2}{-}1$). For $j,k\ins J$, we note $j\succ k$ when $\nu^{(j)}\succ\nu^{(k)}$ (similarly for $\prec$). If $K\sst J$, set
$$K^{\1\succ j}\defs\{k\ins K\mid k\succ j\}\,\,\,\,\h{(similarly for}\,\,\,\,K^{\prec j}).$$
We say that $j\ins K\sst J$ is $\prec$-{\it minimal\1} if $K^{\prec j}\eq\{j\}$. Set
$$V_K\defs\{u\eq(u_j)\ins\C^J\mid u_j\eq0\,\,(j\notins K)\}.$$
Let $V^{(j)}\sst V$ be the subspace on which $r_{\!j+\delf}\eq1$ (and $r_{\!j+\delf}\eq0$ outside $V^{(j)}$). Theorem~\hl{T2}{2} and Corollary~\hl{C1}{1} imply the following.
\par\htt{C2}{}\msn
{\bf Corollary~2.} {\it Assume condition~{\rm(\hl{M1}{M1})} for $\al_{f,k}\notins\Z$ with $\alt_f\gess\tfrac{n}{2}{-}1$.
\skn
{\rm(i)} If $k\eq j\pl\delf$ so that $\om_k\eq[h_j\ddd x]$ $($see {\rm(\hl{6}{6}))}, the shift of $\al_{f,k}\eq\al_{\bf w}([h_j\ddd x])$ in Corollary~{\rm\hl{C1}{1}} depends only on the $u_{j'}$ for $j'\less j$.
\skn
{\rm(ii)} For $j\ins J$, the subspace $V^{(j)}\sst V$ is defined by the equations $u_{j'}\eq p_{j'}$ for $j'\ins J^{\prec j}$. Here the $p_{j'}$ are weighted homogeneous polynomials in variables $u_{j''}$ for $j''\slt j'$ with weighted degree $\ga_{j'}$, and are given by $g_{k,l}^{(1)}\pl u_{j'}$ with $k\eq j'''\pl\delf$, $l\eq j\pl\delf$, and $\nu^{(j')}\pl\nu^{(j''')}\eq\nu^{(j)}$. In particular, the $V^{(j)}$ are smooth with codimension $|J^{\prec j}|$.}
\ms
Note that $V^{(j)}\cap V^{(j')}$ is not necessarily smooth. This implies that the closure of each stratum may have singularities, see \hl{1.5}{1.5} below.
\sk
Let ${\rm SG}(K)\sst\Q_{>0}$ be the {\it semigroup generated additively by\1} the $\ga_k\ins\Q_{>0}$ ($k\ins K$). From Theorem~\hl{T2}{2} we can deduce the following.
\par\htt{P2}{}\msn
{\bf Proposition~2.} {\it Assume condition~{\rm(\hl{M1}{M1})} for $\al_{f,k}\notins\Z$ with $\alt_f\gess\tfrac{n}{2}{-}1$. Let $j\ins K\sst J$.
\skn
{\rm (i)} We have $r_{\!j'+\delf}\eq1$ for any $j'\in J^{\1\succ j}$, if $\ga_j\notins{\rm SG}(K{\setminus}\{j\})$ and $u\ins V_K$ with $u_j\nes 0$.
\skn
{\rm (ii)} We have $r_{\!j'+\delf}\eq1$ for any $j'\in J^{\1\succ j}$, if $\ga_j\ins{\rm SG}(K{\setminus}\{j\})$ and the $u_k$ for $k\ins K{\setminus}\{j\}$ are sufficiently general with $u_{l}$ for $l\notins K{\setminus}\{j\}$ fixed $($even if $j'\notins K$ and $u_{j'}\eq0)$.
\skn
{\rm (iii)} We have $r_{\!j+\delf}\eq0$ for some $u\in V_K$ with $u_j\nes0$ in the case $\ga_j\ins{\rm SG}(K{\setminus}\{j\})$ and $j\ins K$ is $\prec$-minimal. $($Recall that $\delf\eq\mu_f\mi{\rm mod}_{f_1}$, see {\rm(\hl{6}{6}).)}}
\ms
Here ``sufficiently general" means that it is contained in a non-empty Zariski-open subset.
Proposition~\hl{P2}{2} shows the importance of the (partial) additive {\it semigroup structure\1} of the weights $\ga_j$ $(j\ins J)$ for the determination of the $r_{\!j+\delf}$. It gives the (first) {\it affine\1} stratification of the parameter space of the miniversal $\mu$-constant deformation of $f_1$ with coordinates $u_j$ ($j\ins J$). Its strata correspond to {\it bistable\1} subsets $K\sst J$ (where $K$ may be $J$ or $\emptyset$). Here a subset $K\sst J$ is called {\it bistable\1} if the following two conditions are satisfied:
\msn
\vbox{\nin
(a) If $j\ins J$, $k\ins K$, and $j\succ k$, then $j\ins K$.
\skn
(b) If $j\ins J$ with $\ga_j\ins{\rm SG}(K)$, then $j\ins K$.}
\msn
These are called respectively the {\it upper\1} and {\it semigroup\1} stability conditions. The closure of the stratum corresponding to $K$ coincides to $V_K$ (see Corollary~\hl{C3}{3} below), and we delete the closed subspaces corresponding to {\it bistable proper\1} subsets of $K$, that is, the stratum corresponding to $K$ is given by
$$V_K^{\circ}\defs V_K\stm\mcup_{K'\subset K}V_{K'},$$ 
where $K'\sst K$ runs over bistable proper subsets of $K$. Note that bistable subsets are stable by intersections, and the closure of the corresponding stratum is compatible with intersections. 
\par\htt{R1}{}\msn
{\bf Remark~1.} We can determine the bistable subsets $K$ of $J$ by {\it decreasing\1} induction on $|K|$. We first determine those with $|K|\eq|J|{-}1$ by deleting each element from $J$ and verifying the two conditions of bistability, where the obtained bistable subsets are ordered increasingly using the identification $J=\{1,\dots,{\rm mod}_{f_1}\}$ such that the $\ga_j$ are increasing. We say that an element is {\it removable\1} if its complement is bistable. There is at least one removable element, since the two conditions are satisfied by deleting $j$ with $\ga_j$ minimal. We apply the same to the obtained bistable subsets. If we get an already obtained bistable subset, it is of course neglected. We can then proceed by decreasing induction on $|K|$. Note that for any proper bistable subset $K'$ of a bistable subset $K$, there is a removable element of $K$ not contained in $K'$. (Take an element $j\ins K\stm K'$ with $\ga_j$ minimal.) This implies that the strata of the first stratification are affine varieties.
\sk
As a consequence of Theorem~\hl{T2}{2} and Proposition~\hl{P2}{2}, we have the following.
\par\htt{C3}{}\msn
{\bf Corollary~3.} {\it Assume condition~{\rm(\hl{M1}{M1})} for $\al_{f,k}\notins\Z$ with $\alt_f\gess\tfrac{n}{2}{-}1$. If $K\sst J$ is a {\it bistable\1} subset, we have at a sufficiently general point of $V_K$}
\htt{8}{}
$$r_{\!j+\delf}\eq 1\iff j\ins K.
\leqno(8)$$
\ms
On certain locally-closed subspaces of $V_K$, however, the equivalence~(\hl{8}{8}) can hold only after replacing $K$ with a suitable subset of $K$. We thus have to consider a {\it further\1} stratification of each stratum of the first affine stratification, depending on the (partial) additive {\it semigroup structure\1} of the corresponding bistable subset $K\sst J$. In simple cases as in \cite{Kat1}, \cite{Kat2}, \cite{Ca}, where the semigroup structure is not quite complicated, one can easily verify that the BS stratification of the parameter space $V$ of the miniversal $\mu$-constant deformation is described completely by using the bistable subsets of $J$, see \hl{2.1}{2.1}--\hl{2.2}{2} below. In general it does not seem very clear whether each stratum is smooth. (The relation to the stratification by Tjurina numbers seems unclear.)
\sk
We say that a root of a BS polynomial of a weighted homogeneous polynomial $f_1$ with an isolated singularity is {\it shiftable\1} if it is {\it not\1} a root of the BS polynomial of some $\mu$-constant deformation of $f_1$, which is given by a semi-weighted-homogeneous polynomial $f$, see \cite{Va3}. It is well known (and easy to show) that this condition is equivalent to that the root up to sign is strictly greater than $\alt_f{+}1$, assuming condition~(\hl{M1}{M1}) for $\al_{f,k}\notins\Z$ with $\alt_f\gess\tfrac{n}{2}{-}1$, see also \cite{len}. (This does not hold without assuming (\hl{M1}{M1}), for instance if $f_1\eq x^7y\pl xy^5$.) We denote by $R^{\,\rm sh}_{f_1}\,({\cong}\,J)$ the set of {\it shiftable roots\1} up to sign of $b_{f_1}(s)$. (Here $R^{\,\rm sh}_{f_1}\sst\Rt_{f_1}$, since 1 is unshiftable.) 
\sk
Let $\ga_j$ ($j\ins J_{\rm mg}$) be the {\it minimal generators\1} of ${\rm SG}(J)$. The subset $J_{\rm mg}\sst J$ is unique considering the minimal element of ${\rm SG}(J)$ not contained in the subsemigroup
$$\msum_{j\in J_{\rm mg},\,j<k}\,\Z_{>0}\1\ga_j,$$
by induction on $k$. (Note that ${\rm SG}(J)$ is an ordered semigroup such that $\ga_j\pl\ga_{j'}>\ga_j$ for any $j,j'\ins J$.)
\sk
By Theorem~\hl{T2}{2} and Proposition~\hl{P2}{2} we get the following.
\par\htt{C4}{}\msn\vbox{\nin
{\bf Corollary~4.} {\it Assume condition~{\rm(\hl{M1}{M1})} for $\al_{f,k}\notins\Z$ with $\alt_f\gess\tfrac{n}{2}{-}1$. Then all the shiftable roots of $b_f(s)$ are shifted if $u_j\slt0$ $(j\ins J_{\rm mg})$ and $u_j\less 0$ $(j\ins J\stm J_{\rm mg})$. In particular, the shiftable roots are all shifted by adding to $f_1$ only one monomial $h_1$ corresponding to $\ga_1$ in the case all the $\ga_j$ $(j\ins J)$ are contained in the additive semigroup generated by $\ga_1$, for instance, if $\alt_f\pl1\pl\tfrac{1}{{\rm Ord}(T)}$ is a spectral number of $f_1$, where ${\rm Ord}(T)$ is the order of the monodromy $T$.}}
\ms
Here the coefficient of $h_1$ in $f$ can be changed using the $\C^*$-action in case only one monomial is added to $f_1$. To verify the last hypothesis on the spectral number, one can calculate the spectral numbers of $f_1$ using only the weights $w_i$ of variables $x_i$, see (\hl{1.1.20}{1.1.20}) and Remark~\hl{R2.8b}{2.8b} below for such an example.
\sk
We are interested in the following.
\par\htt{Pm1}{}\msn
{\bf Problem 1.} Is it possible that only one shiftable root up to sign $\al$ is unshifted and all the other shiftable ones are shifted, although $\al$ is close to the maximal spectral number?
\ms
This kind of phenomenon is often observed if $\al$ is close to $\alt_f{+}1$, for instance, if the root is associated with a removable element of $J$. In general it can happen in the following case: there is $j_0\ins J$ which is not {\it under\1} $j_{\al}\ins J$ corresponding to $\al$, and if $\ga_j/\ga_{j_0}\ins\Z$ for any $j\ins J$ {\it not over\1} $j_0$. In this case we get the inductive relation that the value of the variable $u_j$ corresponding to each $j\ins J$ {\it under\1} $j_{\al}$ is equal to the value of a certain weighted homogeneous polynomial in variables of {\it strictly lower weights\1} (given by some $g_{k,l}^{(1)}$), see for instance Conjecture~\hl{Cn1}{1} below, where $j_0\eq1$, $\ga_1\eq\tfrac{1}{ab-1}$. Note that all the examples in this paper are computed by using this kind of relation.
\sk
The main problem here is that some of the other shiftable roots could be {\it unshifted.} This actually occurs in the case of deformations of polynomials of BP type, see for instance \hl{2.3}{2.3} below, where the subspace for $\tfrac{83}{56}$ is contained in that for $\tfrac{75}{56}$. This can happen since the action of $\dd_t$ is defined by using only the division by $f_{1,x}\eq ax^{a-1}$ with the action of $\dd_x$ or the one by $f_{1,y}\eq by^{b-1}$ with $\dd_y$. If $f_1$ is {\it nearly BP type\1} which is not BP, this problem does not usually occur, since the action of $\dd_t$ is defined in a more complicated way, although there still remains a certain problem, see Conjecture~\hl{Cn1}{1} below. Note that it is not easy to get an {\it explicit expression\1} of $f$ because of the {\it integer overflow\1} problem especially when $\al$ is close to the maximal spectral number.
\sk
In the case where $\al$ is the unique unshifted shiftable root up to sign of $b_f(s)$, we define ${\rm SR}(f,\al)\defs|R_{f_1}^{\,{\rm sh},<\al}|/|R^{\,\rm sh}_{f_1}|$, called the solitude ratio, where $R_{f_1}^{\,{\rm sh},<\al}=R^{\,\rm sh}_{f_1}\cap(0,\al)$. Since one cannot get this ratio by computing only the BS polynomial of $f$ (here one has to compare it with that of $f_1$), one may also consider the solitude distance ${\rm SD}(f,\al)$ which is the difference between $\al$ and the nearest root up to sign of $b_f(s)$. As these numbers become large, we have more complexity of the defining equations of the subspace on which $\al$ is the unique unshifted shiftable root up to sign of the BS polynomial. We have the following.
\par\htt{Pm2}{}\msn
{\bf Problem 2.} For any $\ep\sgt0$, is there a $\mu$-constant deformation $f$ of a weighted homogeneous polynomial with an isolated singularity which has the unique unshifted shiftable root $\al$ up to sign of $b_f(s)$ with ${\rm SR}(f,\al)\sgt1{-}\ep$? (Similarly for ${\rm SD}(f,\al)$.)
\ms
Examples with ${\rm SR}(f,\al)\eq\tfrac{1}{2}$ are known, see \hl{2.1}{2.1} and \hl{2.4}{2.4} below. We can find an example with ${\rm SR}(f,\al)\eq\tfrac{9}{16}\sgt\tfrac{1}{2}$ for $n\eq3$, see \hl{2.7}{2.7} below. Extending the algorithm to polynomials of {\it nearly\1} BP type, we can get an example with ${\rm SR}(f,\al)\eq\tfrac{2}{3}$ for $n\eq 2$, see \hl{2.6}{2.6} below. This can be extended to the following.
\par\htt{Cn1}{}\msn\vbox{\nin
{\bf Conjecture~1.} Let $f_1\eq x^ay\pl xy^b,$ $\,f\eq f_1\pl x^by^2\pl\msum_{i<b,\,j<b,\,i+j>a}\,u_{i,j}x^iy^j$ with $a\eq b{+}1$. For infinitely many $b\gess5$, the rational number $\al'\defs\tfrac{b(2b-1)}{ab-1}$ is the unique unshifted shiftable root up to sign of $b_f(s)$ choosing the $u_{i,j}\ins\C$ appropriately, where ${\rm SR}(f,\al')\eq1{-}\tfrac{2}{(b-1)(b-2)/2}$ and ${\rm SD}(f,\al')\eq1{-}\tfrac{5b-4}{ab-1}$.}
\ms
The number $\al'$ is the {\it second\1} largest spectral number of $f$, and is associated with $x^{b-1}y^{b-1}$, which is {\it not over\1} $x^by^2$, see (\hl{1.1.16}{1.1.16}) below. It is easy to see that $\al'$ is unshifted if the $u_{i,j}$ coincide with certain rational numbers $c_{i,j}$ which are determined inductively by using certain weighted homogeneous polynomials $h_{i,j}$ of variables with lower weights (defined by a linear function $\ell(i,j)$). Here $h_{i,j}$ is given by $g_{k,l}^{(1)}$ in Theorem~\hl{T2}{2} with monomial of degree 1 deleted (and $k$, $l$ correspond respectively to $x^{b-1-i}y^{b-1-j}$ and $x^{b-1}y^{b-1}$). In this case it is very much expected that the other shiftable roots are shifted. It is, however, rather difficult to prove the last assertion (or to find a counterexample). Indeed, we have to verify that $c_{i,j}$ is different from the value $\ct_{i,j}$ of a weighted homogeneous polynomial $h'_{i,j}$ (which is given by $g_{1,l'}$ with monomial of degree 1 deleted, where $l'$ corresponds to $x^iy^j$) at the $c_{i',j'}$ (but not $\ct_{i'j'}$) for $(i',j')\ins J$ with $\ell(i',j')\slt\ell(i,j)$. (So the $\ct_{i,j}$ are {\it not\1} defined inductively using the $h'_{i,j}$.) The possibility of coincidence seems apparently very slim considering their definitions, but it is quite difficult to show exactly this non-coincidence. One difficulty comes from signs, and this cannot be avoided by replacing $u_j$ with $-u_j$, since signs reappear after {\it substitutions.} (Note that {\it all the shiftable roots are shifted\1} in case every $u_{i,j}$ vanishes, since the $\ga_k$ ($k\ins J$) are contained in the additive semigroup generated by $\ga_1$ and we have the non-vanishing of the coefficient in $f$ of the monomial $x^by^2$, which has the lowest weight in $J$.) There is a similar conjecture for the case $f_1\eq x^a\pl xy^b$ with $a\eq b{+}2$, see also Remark~\hl{R2.6}{2.6} below.
\sk
Note finally that the above argument can be extended to the homogeneous case where condition~(\hl{M1}{M1}) fails. We get an example of a semi-homogeneous polynomial such that the roots of its local BS polynomial are {\it nonconsecutive\1} (that is, not given by the intersection of a connected interval in $\R$ and $\tfrac{1}{d}\1\Z$ with $d$ the degree of the lowest homogeneous part of $f$), for instance,
\vskip-7mm
$$f\eq\tfrac{1}{10}\1x^{10}\pl\tfrac{1}{10}\1y^{10}\pl x^3y^8\pl x^8y^3\pl 8\1x^6y^6\mi 128\1x^7y^7,$$
where $\Rt_f\eq\bl\{\tfrac{2}{10},\dots,\tfrac{14}{10}\br\}\,{\cup}\,\bl\{\tfrac{16}{10}\br\}$ according to Singular \cite{Sing}, see \hl{2.5}{2.5} below. Note that $\Rt_f\sst\tfrac{1}{d}\1\Z\cap(0,n)$ if $f$ is semi-homogeneous and $d$ is the degree of the lowest homogeneous part which has an isolated singularity. (It does not seem easy to find an example as above with $d\less 9$, $n\eq 2$.) This situation is entirely different from the case of homogeneous polynomials with non-isolated singularities, where the roots supported at the origin are {\it consecutive\1} as far as calculated, see \cite{bha}, \cite{nwh}.
\sk
In Section 1 we describe the algorithm after reviewing some basics of Brieskorn lattices. In Section 2 we explain some interesting examples. In Appendix we give some sample codes to compute the $g_{k,l}^{(1)}$ and the bistable subsets.
\sk
This work was partially supported by JSPS Kakenhi 15K04816.
\bs\bs
\vbox{\centerline{\bf 1. Description of the algorithm}
\bsn
In this section we describe the algorithm after reviewing some basics of Brieskorn lattices.}
\par\htt{1.1}{}\msn
{\bf 1.1.~Brieskorn lattices.} Let $f\in\C\{x\}$ be a convergent power series of $n$ variables having an isolated singularity at 0 with $f(0)\eq0$ and $n\gess 2$. Set $(X,0)\defs(\C^n,0)$. The Brieskorn lattice $\Hppf$ (see \cite{Br}) is defined by
\htt{1.1.1}{}
$$\Hppf:=\Om_{X,0}^n/\ddd f{\wedge}\1\ddd\Om_{X,0}^{n-2}.
\leqno(1.1.1)$$
This is a free module of rank $\mu_f$ (with $\mu_f$ the Milnor number of $f$) over $\C\{t\}$ and also over $\Cdti$, see for instance \cite{bl}. We can define the actions of $t$ and $\dti$ respectively by multiplication by $f$ and
\htt{1.1.2}{}
$$\dti[\om]\eq[\ddd f{\wedge}\eta]\q\h{if}\q\ddd\eta\eq\om\q\h{for}\q\om\ins\Om_{X,0}^n,\,\eta\ins\Om_{X,0}^{n-1}.
\leqno(1.1.2)$$
The {\it Gauss-Manin system\1} can be defined by
$$\G_f:=\Hppf[\dd_t],$$
which is the localization of $\Hppf$ by the action of $\dti$. This is a regular holonomic $\D_{X,0}$-module with quasi-unipotent monodromy, and has the $V$-filtration of Kashiwara and Malgrange {\it indexed by\1} $\Q$ so that $\dd_tt\mi\al$ is nilpotent on the graded quotients $\Gr_V^{\al}\G_f$ ($\Q\ins\Q$). The latter is identified with
\htt{1.1.3}{}
$$\G_f^{(\al)}\defs{\rm Ker}\bl((\dd_tt\mi\al)^i\,{:}\,\G_f\tos\G_f\br)\q(i\gg0),
\leqno(1.1.3)$$
and we have the inclusion (see for instance \cite{bl}):
\htt{1.1.4}{}
$$\G_f\into\mprod_{\al\in\Q}\,\G_f^{(\al)}.
\leqno(1.1.4)$$
This implies for $\om\ins\G_f$ the {\it asymptotic expansion\1}
\htt{1.1.5}{}
$$\om=\msum_{\al\in\Q}\,\om^{(\al)}\q\h{with}\q\om^{(\al)}\ins\G_f^{(\al)}\,\,(\al\ins\Q).
\leqno(1.1.5)$$
Set
\htt{1.1.6}{}
$$\al_V(\om):=\min\{\al\in\Q\mid\om^{(\al)}\nes 0\},\q\Gr^V\!\1\om:=\om^{(\al_V(\om))}.
\leqno(1.1.6)$$
\sk
We have the canonical isomorphism
\htt{1.1.7}{}
$$\Gr_V^{\al}\G_f=H^{n-1}(F_{\!f},\C)_{\la}\q(\la\eq e^{-2\pi i\al}),
\leqno(1.1.7)$$
with $F_{\!f}$ the Milnor fiber, and moreover
\htt{1.1.8}{}
$$\aligned&\Gr_V^{\al}\Hppf=F^{n-1-p}H^{n-1}(F_{\!f},\C)_{\la}\\&\h{for}\,\,\,\,\al\eq\be\pl p,\,\be\in(0,1],\,p\ins\Z,\endaligned
\leqno(1.1.8)$$
see \cite[(4.5)]{SS}, \cite{Va1} (and also \cite[(2.6.3)]{bl}).
\par\htt{R1.1a}{}\msn
{\bf Remark 1.1a.} In the {\it weighted homogeneous\1} case, the variable $x_i$ has weight $\om_i$ so that the weighted degree of $f$ is 1, and the filtration $V$ is induced by the filtration on $\Om_{X,0}^n$ by the {\it weighted degree,} where the weight of $\ddd x_i$ is $w_i$. Indeed, we have the Euler field $\xi\eq\msum_{i=1}^n\,w_ix_i\dd_{x_i}$ such that $\xi(f)\eq f$, and
\htt{1.1.9}{}
$$\ddd(\iota_{\xi}\om)=L_{\xi}\1\om,\q\q\ddd f{\wedge}\iota_{\xi}\1\om=f\om\q\q(\om\ins\Om_{X,0}^n),
\leqno(1.1.9)$$
where $L_{\xi}$, $L_{\xi}$ denote respectively the {\it interior product\1} and the {\it Lie derivation\1} respectively, see also \cite[1.1.7]{len}. This calculation implies that the action of $t\dd_t$ on $\G_f$ is semisimple.
\par\htt{R1.1b}{}\msn
{\bf Remark 1.1b.} Assume $f\eq f_1\pl f_{>1}$ is a {\it semi-weighted-homogeneous\1} deformation of a weighted homogeneous polynomial $f_1$ having an isolated singularity at 0. Let $V$ be the decreasing filtration on $\Om_{X,0}^n$ defined by the condition that the weighted degree is at least $\al$. This induces the $V$-filtration on the Gauss-Manin system $\G_f$, see for instance \cite{len}. Moreover we have the canonical isomorphism as graded $\Gr_V^{\ssb}\D_{\C,0}$-modules:
\htt{1.1.10}{}
$$\Gr_V^{\ssb}\G_f=\Gr_V^{\ssb}\G_{f_1}.
\leqno(1.1.10)$$
(This can be shown for instance considering $\Gr_V^{\ssb}\Om_{X,0}^n/\Gr_V^1\ddd f{\wedge}\ddd\Gr_V^{\ssb}\Om_{X,0}^{n-2}$.) This isomorphism makes the calculation of $\Gr_V\dd_t$ simple in the case $f_1$ is BP or nearly BP type with $n\eq2$ as in Remark~\hl{R1.1c}{1.1c} just below.
\par\htt{R1.1c}{}\msn
{\bf Remark~1.1c.} Set $\om^{i,j}\defs x^{i-1}y^{j-1}\ddd x{\wedge}\ddd y\,$ ($i,j\ins\Z_{>0}$) with $n\eq2$.
\sk
For $f_1\eq x^a\pl y^b$, the graded images $\,\dd_t\Gr_V^{\al}[\1\om^{i,j}\1]\eq\Gr_V^{\al-1}\dd_t[\1\om^{i,j}\1]\ins\G_f^{(\al-1)}$ are given by
\htt{1.1.11}{}
$$\dd_t\Gr_V^{\al}[\1\om^{i,j}\1]=\begin{cases}a^{-1}(i{-}a)\1\Gr_V^{\al-1}[\1\om^{i-a,j}\1]&\q(i\sgt a),\\ b^{-1}(j{-}b)\1\Gr_V^{\al-1}[\1\om^{i,j-b}\1]&\q(j\sgt b),\end{cases}
\leqno(1.1.11)$$
with $\al\defs(b\1i{+}aj)/ab$. Here $a^{-1},b^{-1}$ are omitted in the case $f_1\eq\tfrac{1}{a}\1x^a\pl \tfrac{1}{b}\1y^b$.
\sk
In the {\it chain\1} type case, that is, for $f_1\eq x^a\pl xy^b$, we have
\htt{1.1.12}{}
$$\dd_t\Gr_V^{\al}t[\1\om^{i,j}\1]=\begin{cases}a^{-1}(i{-}a{-}\tfrac{j}{b})\1\Gr_V^{\al-1}[\1\om^{i-a,j}\1]&\q(i\sgt a),\\ b^{-1}(j{-}b)\1\Gr_V^{\al-1}[\1\om^{i-1,j-b}\1]&\q(j\sgt b),\end{cases}
\leqno(1.1.12)$$
with $\al\defs(b\1i{+}(a{-}1)j)/ab$. Here $a^{-1},b^{-1}$ are omitted in the case $f_1\eq\tfrac{1}{a}\1x^a\pl \tfrac{1}{b}\1xy^b$.
\sk
In the {\it loop\1} type case, that is, for $f_1\eq x^ay\pl xy^b$, we have
\htt{1.1.13}{}
$$\dd_t\Gr_V^{\al}[\1\om^{i,j}\1]=\begin{cases}(ab{-}1)^{-1}(b(i{-}a){-}j{+}1)\1\Gr_V^{\al-1}[\1\om^{i-a,j-1}\1]&\q(i\sgt a),\\(ab{-}1)^{-1}(a(j{-}b){-}i{+}1)\1\Gr_V^{\al-1}[\1\om^{i-1,j-b}\1]&\q(j\sgt b),\end{cases}
\leqno(1.1.13)$$
with $\al\defs((b{-}1)i{+}(a{-}1)j)/(ab{-}1)$, see (\hl{1.1.14}{1.1.4}) just below. They are multiplied respectively by $a$ and $b$ in the case $f_1\eq\tfrac{1}{a}\1x^ay\pl \tfrac{1}{b}\1xy^b$. Here we assume that $(i,j)$ is a linear combination of $(a,1)$ and $(1,b)$ with positive coefficients, and moreover this holds after replacing $(i,j)$ with $(i{-}a,j{-}1)$ or $(i{-}1,j{-}a)$ so that the coefficients in (\hl{1.1.13}{1.1.13}) are {\it non-negative.} We assume a similar condition for (\hl{1.1.11}{1.1.11}-\hl{1.1.12}{12}).
\par\htt{R1.1d}{}\msn
{\bf Remark~1.1d.} Assume $f$ is a polynomial of nearly BP loop type in two variables, that is, $f\eq x^ay\pl xy^b$. Let $\ell$ be the weighted degree function on $\Z^2$ such that $\ell(a,1)\eq\ell(1,b)\eq1$. It is easy to see that $\ell$ is defined by
\htt{1.1.14}{}
$$\ell(i,j)=(b'i\pl a'j)/(ab{-}1)\q(a\eq a'{+}1,\,\,b\eq b'{+}1).
\leqno(1.1.14)$$
We have $\mu_f\eq ab$, see \cite[A.1]{wh}. The Jacobian ring is spanned over $\C$ by the monomials $x^iy^j$ for $i\slt a$, $j\slt b$. (In the chain type case, that is, for $f_1\eq x^a\pl xy^b$, the Jacobian ring is spanned  by the monomials $x^iy^j$ for $i\slt a$, $j\slt b{-}1$ and $y^{b-1}$, where $\mu_f\eq a(b{-}1)+1$.)
\sk
In the case $a\eq b{+}1$, one can verify that
\htt{1.1.15}{}
$$\dim H^1(F_{\!f},\C)_{\la}\eq\begin{cases}1&(\la^{ab-1}\eq1,\,\la\nes 1),\\ 2&(\la\eq 1),\\ 0&(\la^{ab-1}\nes1).\end{cases}
\leqno(1.1.15)$$
with $\ga_1\eq1/(ab{-}1)$, where the minimal, maximal, and second largest spectral numbers are given respectively by
\htt{1.1.16}{}
$$(2b{-}1)/(ab{-}1),\q(2b^2{-}1)/(ab{-}1),\q(2b^2{-}b)/(ab{-}1).
\leqno(1.1.16)$$
\par\htt{R1.1e}{}\msn
{\bf Remark 1.1e.} In general, a polynomial $f$ is called a polynomial of {\it nearly BP loop type\1} if
$$f\eq\msum_{i=1}^n\,c_ix_i^{a_i}x_{i+1}\q(c_i\ins\C^*,\,a_i\ins\Z_{\ges 2}),$$
where $\{1,\dots,n\}$ is identified with $\Z/n\Z$. (It is called a non-degenerate invertible polynomial or potential of loop type in mirror symmetry, see for instance \cite{Kre}, \cite{EG}.)
\sk
For $n\eq3$, we have $\mu_f\eq abc$ and $T^{abc+1}\eq1$ with $T$ the monodromy, where the $a_i$ are denoted by $a,b,c$. The weights multiplied by $abc{+}1$ are given by
$$(b{-}1)c{+}1,\q(c{-}1)a{+}1,\q(a{-}1)b{+}1.$$
In case their greatest common divisor is 1, an assertion similar to (\hl{1.1.15}{1.1.15}) seems to hold replacing $ab{-}1$ by $abc{+}1$ and 2 by 0 for $\la\eq1$.
\sk
Indeed, for $n\gess2$, let $v^{(i)}\ins\N^n$ with $f_1\defs\msum_{i=1}^n\,x^{v^{(i)}}$ a polynomial of nearly BP loop type. It is easy to see that the $(n{-}1)$-th exterior product
$$\xi^{(k)}\defs\h{$\bigwedge$}_{i\ne k}\1(v^{(i)}{-}\1v^{(k)})\ins\h{$\bigwedge$}^{n-1}\Z^n\cong\Z^n$$
is independent of $k\ins[1,n]$ up to sign. It is a primitive vector, that is, the greatest common divisor of its components is 1, if and only if the parallelotope spanned by the $v^{(i)}\mi v^{(k)}$ ($i\nes k$) has no lattice point except for the vertices, since the condition is equivalent to the existence of $v\ins\Z^n$ such that $v{\wedge}\1\xi^{(k)}$ generates $\bigwedge^{\!n}\Z^n$, that is, the determinant of $\bl(v,v^{(i)}\mi v^{(k)}\,(i\nes k)\br)$ is equal to $\pm1$. This seems to imply condition~(\hl{M1'}{M1})$'$ for $\la\nes1$ using some projection to $\Z^{n-1}$. (This may be known to some specialist.) Note that the components of $\xi^{(k)}$ are given up to sign by the above numbers and the modified exponents (see \hl{A3}{A.3} in Appendix) in the case $n\eq3$ and $2$ respectively.
\par\htt{R1.1f}{}\msn
{\bf Remark 1.1f.} Let $\al_{f,1},\dots,\al_{f,\mu_f}$ be the {\it spectral numbers\1} of $f$, see \cite{St2} (and also \cite{DS}, \cite{JKSY}). They are assumed to be weakly increasing. It is well known (see for instance \cite{SS}, \cite{bl}, \cite{Va1} and also (\hl{1.1.8}{1.1.8}) that we have the equality
\htt{1.1.17}{}
$$\dim_{\C}\Gr_V^{\al}\Om_f^n=\#\{k\ins[1,\mu_f]\mid\al_{f,k}\eq\al\}\q(\forall\,\al\ins\Q),
\leqno(1.1.17)$$
with
\htt{1.1.18}{}
$$\Om_f^n:=\Om_{X,0}^n/\ddd f{\wedge}\Om_{X,0}^{n-1}=\Hppf\!/\dti\Hppf.
\leqno(1.1.18)$$
We have the symmetry of spectral numbers as is well known (see \cite{St2}):
\htt{1.1.19}{}
$$\al_{f,k}\pl\al_{f,l}\eq n\q(k\pl l\eq\mu_f\pl1).
\leqno(1.1.19)$$
\sk
In the case $f$ is a weighted homogeneous polynomial of weights $w_i$, the spectral numbers can be computed by
\htt{1.1.20}{}
$$\msum_{k=1}^{\mu_f}\,t^{\al_{f,k}}=\mprod_{i=1}^n\,(t^{w_i}\mi t)/(1\mi t^{w_i}),
\leqno(1.1.20)$$
see \cite{St2}, \cite[Section 1.5]{JKSY}.
\par\htt{R1.1g}{}\msn
{\bf Remark 1.1g.} By \cite{Ma}, the {\it reduced\1} BS polynomial $b_f(s)/(s{+}1)$ is equal to the minimal polynomial of the action of $-\dd_tt$ on
$$\Htppf\!/t\Htppf.$$
Here $t\Htppf$ may be replaced by $\dti\Htppf$, see for instance \cite{len}.
\par\htt{1.2}{}\msn
{\bf 1.2.~Proof of Theorems~\hl{T1}{1} and \hl{T2}{2}.} We have $f\eq f_1\pl\msum_{j\in J}\,u_jh_j$ as in the introduction, where the $u_j$ are viewed as constants. For $\nu\in\Z_{>0}^n$, set $\om^{\nu}\defs[x^{\nu-\ob}\ddd x]\ins\Hppf$. We can easily verify that
\htt{1.2.1}{}
$$\bl(\dd_tt\mi\al_{\bf w}(\nu)\br)\om^{\nu}=-\msum_{j\in J}\,\ga_ju_j\1\dd_t\1\om^{\nu+\nu^{(j)}},
\leqno(1.2.1)$$
using (\hl{1.1.9}{1.1.9}), where $h_j\eq x^{\nu^{(j)}}$, see also \cite{len}. Comparing the {\it asymptotic expansions\1} of both sides (using (\hl{1.1.10}{1.1.10})), we can determine the theoretically lowest term of the asymptotic expansion of $\om^{\nu}$ by decreasing induction on $\al_{\bf w}(\nu)$. Here what is important is the ratio of the weight of the variable in which one is interested and the minimal weight of the variables, since this gives the maximal number of procedures to which one applies (\hl{1.2.1}{1.2.1}) inductively, see also the sample codes in \hl{A1}{A.1}-\hl{A2}{2} of Appendix for more details. The last assertion of Theorem~\hl{T2}{2} also follows from the above argument. This finishes the proof of Theorems~\hl{T1}{1} and \hl{T2}{2}.
\par\htt{R1.2}{}\msn
{\bf Remark~1.2.} To calculate the weighted homogeneous polynomials $g_{k,l}^{(1)}$ with weighted degree given in Theorem~\hl{T2}{2} in the two variable case, we can apply the iteration of (\hl{1.2.1}{1.2.1}), where we get sequences $\{j_i\}\ins J^r$ with $\msum_{i=1}^r\,\ga_{j_i}$ the weighted degree of $g_{k,l}^{(1)}$. The contribution of this sequence to the coefficient of $\mprod_{i=1}^r\,u_{j_i}$ in $g_{k,l}^{(1)}$ is given by
\htt{1.2.2}{}
$$(-1)^r\,\mprod_{i=1}^r\,\bl(\msum_{j=i}^r\,\ga_j\br)^{-1}\,\mprod_{i=1}^r\,\ga_i,
\leqno(1.2.2)$$
multiplied by
\htt{1.2.3}{}
$$\mprod_{j=1}^{[(p+1)/a]}(p{+}1{-}ja)/a\,\,\mprod_{j=1}^{[(q+1)/b]}(q{+}1{-}jb)/b,
\leqno(1.2.3)$$
(which vanishes if $p{+}1\ins a\1\Z$ or $q{+}1\ins b\1\Z$), where $x^py^q\eq g\1\mprod_{i=1}^r\,h_{j_i}$ if $\om_k=g\1\ddd x{\wedge}\ddd y$. (Here $p,q$ {\it cannot\1} be determined by $\om_k,\om_{l}$ in general.) Note that $\msum_{j=i}^r\,\ga_j$ comes from the action of the differential operator on the left-hand side of (\hl{1.2.1}{1.2.1}), and (\hl{1.2.3}{1.2.3}) is obtained by calculating the action of $\dd_t$, where (\hl{1.1.10}{1.1.10}) is used in an essential way. The divisions by $a$ and $b$ must be omitted if we consider $f_1\eq\tfrac{1}{a}\1x^a\pl\tfrac{1}{b}\1y^b$ instead of $x^a\pl y^b$. Here we have to consider all the possible sequences with various lengths, and we usually get a {\it combinatorial problem.} This formula is used in \cite{len}.
\par\htt{1.3}{}\msn
{\bf 1.3.~Proof of Proposition~\hl{P1}{1}.} The maximal spectral number $\al_{f,\mu_f}$ is equal to $n\mi\alt_f$ by the symmetry of spectral numbers, see (\hl{1.1.19}{1.1.19}). The assertion then follows from the generalized Jordan decomposition, see for instance \cite[Rem.\,A.7c]{rh}. This finishes the proof of Proposition~\hl{P1}{1}.
\par\htt{R1.3}{}\msn
{\bf Remark~1.3.} Let $f\eq\msum_{\be\ges1}\,f_{\be}$ be a semi-weighted-homogeneous polynomial, where the $f_{\be}$ have weighted degree $\be$ and $f_1$ has an isolated singularity at 0. Let $m\ins\Z_{>0}$ such that $m\be\ins\Z$ if $f_{\be}\nes 0$. Then we have the one-parameter family $g_u\defs\msum_{\be\ges1}\,u^{m(\be-1)}f_{\be}$ such that $g_0\eq f_1$, $g_1\eq f$, and $b_{g_u}(s)$ is constant for $u\nes0$.
\par\htt{1.4}{}\msn
{\bf 1.4.~Proof of Proposition~\hl{P2}{2}.} This follows by applying (\hl{1.2.1}{1.2.1}) inductively to the asymptotic expansion of $[x^{\nut}\ddd x]$ with $\nut\eq\nu^{(j')}\mi\nu^{(j)}$. Recall that the sign of the coefficients of $g_{k,l}^{(1)}$ depends only on the total degree of the $u_j$, and we have the positivity if the $u_j$ are replaced by $-u_j$. Here we use (\hl{1.1.10}{1.1.10}) in an essential way. This finishes the proof of Proposition~\hl{P2}{2}.
\par\htt{R1.4}{}\msn
{\bf Remark~1.4.} The minimal root of $b_{f_1}(s)$ is not a root of $b_f(s)$ if $u_j\nes 0$ for some $j\ins J$. Indeed, let $j_0\ins J$ with $u_{j_0}\nes0$ and $\ga_{j_0}$ minimal. Let $h_{j_0}$ be the corresponding monomial. Set $g\defs\mprod_{i=1}^n\,x_i^{e_i-2}/h_{j_0}$. Consider the asymptotic expansion of $[g\ddd x]$, and apply (\hl{1.2.1}{1.2.1})
\par\htt{1.5}{}\msn
{\bf 1.5~Singularities of the closures of BS strata.} In the notation of the introduction, the closure of each stratum of the BS stratification of $V$ is not necessarily smooth. For instance, let $f_1\eq x^9y\pl xy^8$. With the notation of Corollary~\hl{C2}{2}, the intersection $V^{(10)}\cap V^{(11)}$ has non-isolated singularities whose transversal slice is a surface singularity of type $A_4$. This can be verified by typing BC9~8~11 for ./a.out of the code in \hl{A2}{A.2} of Appendix and examining the difference between [2;5,6] and [3;4,7] in the output after the substitution caused by [3;5,6], where Corollary~\hl{C2}{2}\,(ii) is used.
\bs\bs
\vbox{\centerline{\bf 2. Examples}
\bsn
In this section we explain some interesting examples.}
\par\htt{2.1}{}\msn
{\bf 2.1.~Example I.} Let $f_1\eq x^7{+}\1y^5$ or $x^9{+}\1y^4$, see \cite{Kat1}, \cite{Kat2}. The weights $\ga_j$ ($j\ins J$) multiplied by 35 or 36 and the exponents of the corresponding monomials are as follows:
\htt{2.1.1}{}
$$\scalebox{0.9}{$\begin{array}{ccccccc}1&6&11&&\\&&4\end{array}\begin{array}{ccccccc}3,3&4,3&5,3&&&\\&&5,2\end{array}\h{or}\begin{array}{cccccccc}&&&2&6&10&&\\&&&&&1\end{array}\begin{array}{ccccccc}5,2&6,2&7,2\\&&7,1\end{array}$}
\leqno(2.1.1)$$
We see that the number of strata of the first affine stratification is 6 or 5, where the nonempty bistable subsets are as below:
\htt{2.1.2}{}
$$\scalebox{0.8}{$\begin{array}{cccccccc}1&6&11&&\\&&4\end{array}
\begin{array}{cccccccc}&6&11&&\\&&4\end{array}
\begin{array}{cccccccc}&6&11&&\\&&\end{array}
\begin{array}{cccccccc}&&11&&\\&&4\end{array}
\begin{array}{cccccccc}&&11&&\\&&\end{array}$}
\leqno(2.1.2)$$
or
\htt{2.1.3}{}
$$\scalebox{0.8}{$\begin{array}{ccccccccc}2&6&10&&&\\&&1\end{array}
\begin{array}{cccccccc}2&6&10&&&\\&&\end{array}
\begin{array}{cccccccc}&6&10&&&\\&&\end{array}
\begin{array}{cccccccc}&&10&&&\\&&\end{array}$}
\leqno(2.1.3)$$
The open stratum of the first affine stratification corresponding to $K\eq J$ contains one or two subspaces on which the assertion~(\hl{8}{8}) holds with $J$ replaced by $J\stm\{2\}$ or $J\stm\{2\}$ and $J\stm\{3\}$, that is,
\htt{2.1.4}{}
$$\scalebox{0.8}{$
\begin{array}{ccccccc}1&6&11&&\\&&\end{array}\h{or}
\begin{array}{ccccccccccc}&&&6&10&&\\&&&&1\end{array}\h{and}
\begin{array}{ccccccccccc}&&2&\h{\,\,\,}&10&&&\\&&&&1\end{array}$}
\leqno(2.1.4)$$
using Proposition~\hl{P2}{2}. The corresponding subspaces are as follows:
\htt{2.1.5}{}
$$\{u_2\eq c\1u_1^4\}\q\h{or}\q\{u_2\eq c'u_1^2\}\,\,\,\h{and}\,\,\,\{u_2\eq c''u_1^2,\,u_3\eq c'''u_1^6\}.
\leqno(2.1.5)$$
These follow from Theorem~\hl{T2}{2} looking at the weights $\ga_j$ and (\hl{2.1.1}{2.1.1}). Note that $\ga_2\eq 4\ga_1$ or $\ga_2\eq 2\ga_1$ and $\ga_3\eq 6\ga_1$. (We have to add the minimal spectral number $\tfrac{12}{35}$ or $\tfrac{13}{36}$ plus one to the weight $\ga_j$ in order to get the spectral number $\al_{f,k}$ corresponding to $\ga_j$, see (\hl{6}{6}).) It is not necessarily easy to determine the constants $c,c',c'',c'''\in\C^*$. Using a computer, we can get that
\htt{2.1.6}{}
$$c\eq\tfrac{6}{175},\q c'\eq\tfrac{1}{3},\q c''\eq\tfrac{7}{18},\q c'''\eq\tfrac{2429}{1259712}.
\leqno(2.1.6)$$
To determine $c'''$, we use the vanishing of $g_{2,23}^{(1)}\eq{-}u_2\pl\tfrac{7}{18}u_1^2$ (which gives $c''$) and that of
\htt{2.1.7}{}
$$g^{(1)}_{1,23}\eq{-}u_3\mi\tfrac{7}{72}u_2^3\pl\tfrac{7}{27}u_1^2u_2^2\mi\tfrac{175}{1458}u_1^4u_2\pl\tfrac{595}{39366}u_1^6,
\leqno(2.1.7)$$
where $23\eq\mu_f\mi1$. Notice that the {\it sign\1} depends only on the degrees of monomials. This holds in general applying repeatedly (\hl{1.2.1}{1.2.1}). For $c'$, we need the vanishing of $g^{(1)}_{1,22}\eq{-}u_2\pl\tfrac{1}{3}u_1^2$. (For $c$ with $f_1\eq x^7{+}\1y^5$, we use the vanishing of $g^{(1)}_{1,22}\eq{-}u_2\pl\tfrac{6}{175}u_1^4$ with $22\eq\mu_f\mi2$.)
\sk
For a $\mu$-constant deformation $f$ of $f_1\eq x^9{+}y^4$ having the unique unshifted shiftable root $\tfrac{55}{36}$ up to sign of $\bt_f(s)$, the distribution of roots up to sign is as below:
$$\setlength{\unitlength}{0.3cm}
\begin{picture}(50,2.5)
\multiput(40,1.5)(4,0){1}{\circle{.5}}
\multiput(41,1.5)(4,0){3}{\circle{.5}}
\multiput(3,1.5)(4,0){8}{\circle*{.5}}
\multiput(12,1.5)(4,0){7}{\circle*{.5}}
\multiput(21,1.5)(4,0){5}{\circle*{.5}}
\put(45,1.5){\circle*{.5}}
\put(13,1.5){\circle*{.5}}
\put(5,1.5){\circle*{.5}}
\put(4,1.5){\circle*{.5}}
\put(2.4,0){$\scriptscriptstyle\frac{13}{36}$}
\put(36.4,0){$\scriptscriptstyle\frac{47}{36}$}
\put(44.4,0){$\scriptscriptstyle\frac{55}{36}$}
\put(48.4,0){$\scriptscriptstyle\frac{59}{36}$}
\put(39.4,0){$\scriptscriptstyle\frac{50}{36}$}
\end{picture}$$
Here the black and white vertices represent respectively the roots of $b_f(s)$ up to sign and those of $b_{f_1}(s)$ which are not roots of $b_f(s)$ up to sign (that is, shifted). We have ${\rm SR}(f,\tfrac{55}{36})\eq\tfrac{1}{2}$ and ${\rm SD}(f,\tfrac{55}{36})\eq\tfrac{2}{9}$.
\par\htt{2.2}{}\msn
{\bf 2.2.~Example II.} Let $f_1\eq x^7{+}\1y^6$. The weights of parameters multiplied by 42 and the exponents of the corresponding monomials are as follows:
\htt{2.2.1}{}
$$\scalebox{0.85}{$\begin{array}{cccccc}4&10&16\\&3&9&&\\&&2\end{array}
\begin{array}{ccc}3,4&4,4&5,4\\&4,3&5,3\\&&5,2\end{array}$}
\leqno(2.2.1)$$
Following Remark~\hl{R1}{1} in the introduction, we get ten nonempty bistable subsets of $J$ as below:
\htt{2.2.2}{}
$$\scalebox{0.8}{$
\begin{array}{cccccccc}4&10&16&&\\&3&9\\&&2\\ \\&10&16&&&\\&&9\end{array}
\begin{array}{cccccccc}4&10&16&&\\&3&9\\&\\ \\4&10&16&&&\\&\end{array}
\begin{array}{cccccccc}4&10&16&&\\&&9\\&&2\\ \\&10&16&&&&\\&\end{array}
\begin{array}{cccccccc}4&10&16&&\\&&9\\&\\ \\&&16&&&&\\&&9\end{array}
\begin{array}{cccccccc}&10&16\\&3&9\\&\\ \\&&16\\&\end{array}$}
\leqno(2.2.2)$$
The first and third subsets have respectively one and two subsets as below, which are not bistable, but correspond to the sets of {\it shifted\1} roots up to sign of BS polynomial $b_f(s)$ by adding 13 and dividing it by 42:
\htt{2.2.3}{}
$$\scalebox{0.8}{$
\begin{array}{cccccccc}&10&16&&&\\&3&9\\&&2\end{array}\h{and}
\begin{array}{cccccccc}&&&10&16&&&\\&&&&9\\&&&&2\end{array}
\begin{array}{cccccccc}4&&&16&&\\&&&9\\&&&2\end{array}$}
\leqno(2.2.3)$$
These are determined by using the (partial) {\it semigroup structure\1} of $K$ (more precisely, $4\eq2\scd2$ and $10\eq5\scd2\eq4\pl3\scd2\eq2\scd4\pl2$). The corresponding subspaces are respectively as follows:
\htt{2.2.4}{}
$$\{u_3\eq c\1u_1^2\}\,\,\,\h{and}\,\,\,\{u_2\eq0,\,\,u_3\eq c'u_1^2\},\,\,\{u_2\eq0,\,\,u_3\eq c''u_1^2,\,u_5\eq c'''u_1^5\}.
\leqno(2.2.4)$$
The first subset of (\hl{2.2.2}{2.2.2}) (that is, $J$) has only one subset, since $5\ins J$ is {\it over\1} $2\ins J$ with $\ga_5\eq\tfrac{10}{42}$ and $\ga_2\eq\tfrac{3}{42}\notins{\rm SG}(J{\setminus}\{2\})$. We have
\htt{2.2.5}{}
$$c\eq c'\eq\tfrac{2}{7},\q c''\eq\tfrac{5}{14},\q c'''\eq{-}\tfrac{5}{16464}.
\leqno(2.2.5)$$
It is not easy to determine $c'''$. (There is another method used in \cite{len}, see Remark~\hl{R1.2}{1.2}.) In this paper we calculate it by combining the vanishing of $g^{(1)}_{1,38}\eq{-}u_2$, $g^{(1)}_{2,41}\eq{-}u_3\pl\tfrac{5}{14}u_1^2$ (which gives $c''$), and that of
\htt{2.2.6}{}
$$g^{(1)}_{1,41}\eq{-}u_5\mi\tfrac{25}{84}\1u_2^2u_3-\tfrac{25}{84}\1u_1u_3^2+\tfrac{25}{98}\1u_1^2u_2^2+\tfrac{25}{147}\1u_1^3u_3-\tfrac{95}{4116}\1u_1^5,
\leqno(2.2.6)$$
see \hl{A1}{A.1}-\hl{A2}{2} in Appendix. Here $41\eq\mu_f\mi1$. For $c,c'$, we use the vanishing of $g^{(1)}_{1,39}\eq{-}u_3\pl\tfrac{2}{7}u_1^2$. 
\par\htt{2.3}{}\msn
{\bf 2.3.~Example III.} Let $f_1\eq x^8{+}\1y^7$. The weights of parameters multiplied by 56 and the exponents of the corresponding monomials are as follows:
\htt{2.3.1}{}
$$\scalebox{0.85}{$\begin{array}{cccccccccc}5&12&19&26&&\\&4&11&18\\&&3&10&&\\&&&2\end{array}
\begin{array}{ccccccc}3,5&4,5&5,5&6,5\\&4,4&5,4&6,4\\&&5,3&6,3\\&&&6,2\end{array}$}
\leqno(2.3.1)$$
Following Remark~\hl{R1}{1}, we get 24 nonempty proper bistable subsets of $J$ as below:
$$\scalebox{0.8}{$
\begin{array}{cccccccc}5&12&19&26&\\&4&11&18\\&&3&10\\&&&
\\ \\5&12&19&26&\\&&11&18\\&&&10\\&&&
\\&12&19&26&\\&&11&18\\&&&\\&&&\\&&19&26&\\&&&18\\&&&\end{array}
\begin{array}{cccccccc}5&12&19&26&\\&4&11&18\\&&&10\\&&&2
\\ \\&12&19&26&\\&4&11&18\\&&&10\\&&&
\\&12&19&26&\\&&&18\\&&&10\\&&&\\&12&19&26&\\&&&\\&&&\end{array}
\begin{array}{cccccccc}5&12&19&26&\\&4&11&18\\&&&10\\&&&
\\ \\&12&19&26&\\&&11&18\\&&3&10\\&&&
\\&&19&26&\\&&11&18\\&&&10\\&&&\\&&&26&\\&&&18\\&&&10\end{array}
\begin{array}{cccccccc}5&12&19&26&\\&&11&18\\&&3&10\\&&&
\\ \\&12&19&26&\\&&11&18\\&&&10\\&&&
\\&12&19&26&\\&&&18\\&&&\\&&&\\&&19&26&\\&&&\\&&&\end{array}
\begin{array}{cccccccc}&12&19&26&\\&4&11&18\\&&3&10\\&&&
\\ \\5&12&19&26&\\&&&18\\&&&10\\&&&
\\&&19&26&\\&&11&18\\&&&\\&&&\\&&&26&\\&&&18\\&&&\end{array}
\begin{array}{cccccccc}&12&19&26&\\&4&11&18\\&&&10\\&&&2
\\ \\&12&19&26&\\&4&11&18\\&&&\\&&&
\\&&19&26&\\&&&18\\&&&10\\&&&\\&&&26&\\&&&\\&&&\end{array}$}$$
It is not easy to determine the {\it finer\1} stratification of the first affine stratification. We can verify that three roots of $\bt_f(s)$ corresponding to the largest three spectral numbers $\al_{f,42}\eq\tfrac{97}{56}$, $\al_{f,41}\eq\tfrac{90}{56}$, $\al_{f,40}\eq\tfrac{89}{56}$ have no proper non-empty closed subspaces of affine strata on which some of these roots are {\it unshifted.} Here $\al_{f,k}$ corresponds to $j\ins J\eq[1,10]$ if $k\mi j\eq32\,({=}\,\mu_f{-}10)$. (For instance, if the root $\tfrac{90}{56}$ is {\it unshifted,} then we must have $u_2\eq0$, $u_4\eq0$, $u_6\eq 0$ inductively using Proposition~\hl{P2}{2}, where $\ga_2\eq\tfrac{3}{56}$, $\ga_4\eq\tfrac{5}{56}$, $\ga_6\eq\tfrac{11}{56}$. The other weights smaller than $\tfrac{19}{56}$ are {\it even\1} integers divided by 56, and $\tfrac{19}{56}$ is not contained in the semigroup generated by them. We thus get that $u_9\eq0$ with $9\eq41{-}32$. See Remark~\hl{R1.4}{1.4} for $\tfrac{97}{56}$.) These make the computer calculation quite simple allowing us to avoid {\it integer overflow.} We then see that the roots between $\tfrac{75}{56}$ and $\tfrac{83}{56}$ up to sign (corresponding to $j\ins[3,7]$) are {\it unshifted\1} on some non-empty proper closed subspaces of certain affine strata. For instance, the root $\tfrac{83}{56}$ up to sign is {\it unshifted\1} on the subspace defined by the vanishing of
$$\aligned g_{3,39}^{(1)}&\eq{-}u_3\pl\tfrac{5}{16}u_1^2,\q g_{2,39}^{(1)}\eq{-}u_4\pl\tfrac{5}{8}u_1u_2,\\ g_{1,39}^{(1)}&\eq{-}u_7\pl\tfrac{5}{8}u_1u_5\mi\tfrac{5}{56}u_3^3\mi\tfrac{15}{28}u_2u_3u_4\mi\tfrac{15}{56}u_1u_4^2\pl\tfrac{65}{1792}u_2^4\pl\tfrac{195}{448}u_1u_2^2u_3\\&\q\pl\tfrac{195}{896}u_1^2u_3^2\pl\tfrac{195}{448}u_1^2u_2u_4\mi\tfrac{195}{1024}u_1^3u_2^2\mi\tfrac{195}{2048}u_1^4u_3\pl\tfrac{377}{32768}u_1^6,\endaligned$$
(see \hl{A1}{A.1}-\hl{A2}{2} in Appendix), that is, on the subspace
$$\bl\{u_3\eq\tfrac{5}{16}u_1^2,\,\,u_4\eq\tfrac{5}{8}u_1u_2,\,\,u_7=\tfrac{1}{3584}u_1^6\pl\tfrac{15}{1792}u_1^3u_2^2\pl\tfrac{65}{1792}u_2^4\pl\tfrac{5}{8}u_1u_5\br\}.$$
Similarly the roots $\tfrac{82}{56}$, $\tfrac{81}{56}$, $\tfrac{76}{56}$, $\tfrac{75}{56}$ up to sign are {\it unshifted\1} respectively on the subspaces
$$\aligned&\bl\{u_2\eq0,\,\,u_3=\tfrac{3}{8}u_1^2,\,\,u_6\eq{-}\tfrac{5}{112}u_1^3u_4\br\},\q\bl\{u_1\eq u_2\eq0,\,\,u_5\eq\tfrac{2}{7}u_4^2\br\},\\&\bl\{u_4\eq\tfrac{1}{2}u_1u_2\br\},\q\bl\{u_3=\tfrac{5}{16}u_1^2\br\}.\endaligned$$
Indeed, $g^{(1)}_{1,38}$ and $g^{(1)}_{1,37}$ are given respectively by
$$-u_6\mi\tfrac{15}{56}u_2u_3^2\mi\tfrac{15}{56}u_2^2u_4\mi\tfrac{15}{28}u_1u_3u_4\pl\tfrac{5}{32}u_1u_2^3\pl\tfrac{15}{32}u_1^2u_2u_3\pl\tfrac{5}{32}u_1^3u_4\mi\tfrac{55}{512}u_1^4u_2,$$
$$-u_5+\tfrac{2}{7}u_4^2\mi\tfrac{1}{4}u_2^2u_3\mi\tfrac{1}{4}u_1u_3^2\mi\tfrac{1}{2}u_1u_2u_4\pl\tfrac{15}{64}u_1^2u_2^2\pl\tfrac{5}{32}u_1^3u_3\mi\tfrac{23}{1024}u_1^5.$$
Note that $\tfrac{83}{56},\tfrac{82}{56},\tfrac{81}{56}$ {\it cannot\1} be the unique unshifted shiftable root up to sign. (Indeed, in the latter two cases, the subspace is contained in a coordinate hyperplane. In the first case, it is contained in the subspace for $\tfrac{75}{56}$.)
\par\htt{R2.3}{}\msn
{\bf Remark~2.3.} One can examine the above computation using Singular \cite{Sing} as follows.
\ms
\vbox{\tiny\sf\pv@LIB "gmssing.lib"; ring R=0,(x,y),ds; poly a=2/3;@
\pv@poly u_1=a; poly u_2=a; poly u_3=5/16*u_1^2; poly u_4=5/8*u_1*u_2;@
\pv@poly u_5=a; poly u_6=a; poly u_8=a; poly u_9=a; poly u_10=a;@
\pv@poly u_7=1/3584*u_1^6+15/1792*u_1^3*u_2^2+65/1792*u_2^4+5/8*u_1*u_5;@
\pv@poly f=x^8+y^7+u_1*x^6*y^2+u_2*x^5*y^3+u_3*x^4*y^4+u_4*x^3*y^5+@
\pv@u_5*x^6*y^3+u_6*x^5*y^4+u_7*x^4*y^5+u_8*x^6*y^4+u_9*x^5*y^5+u_10*x^6*y^5;@
\pv@bernstein(f);@}
\skn
Here the $u_j$ for $j\nes3,4,7$ can be arbitrary rational numbers (as long as they are not too much complicated for Singular). One should always get a root $\tfrac{83}{56}$ up to sign together with $\tfrac{75}{56}$.
\par\htt{2.4}{}\msn
{\bf 2.4.~Example IV.} Let $f_1\eq x^9\pl y^7$. The weights of parameters multiplied by 63 and the exponents of the corresponding monomials are as follows:
$$\scalebox{0.85}{$\begin{array}{cccccccccc}3&10&17&24&31&&\\&1&8&15&22\\&&&6&13\\&&&&4\end{array}\begin{array}{cccccccccc}3,5&4,5&5,5&6,5&7,5\\&4,4&5,4&6,4&7,4\\&&&6,3&7,3\\
&&&&7,2\end{array}$}$$
There are 34 nonempty bistable subsets $K$ of $J$, where the numbers of $K$ with $|K|\eq i$ are $1, 1, 2, 3, 4, 4, 5, 4, 4, 3, 2, 1$ for $i\eq12,\dots,1$ respectively. If $u_1\nes0$, the roots corresponding to $j\ins J$ {\it over\1} $1$ are all shifted by Proposition~\hl{P2}{2}, and the calculation is not very difficult when $u_1\eq0$. So we examine the shift of the root up to sign $\tfrac{92}{63}$ corresponding to $j\eq 7$. We have $\ga_1\eq\tfrac{1}{63}$, $\ga_7\eq\tfrac{13}{63}$, and their ratio is $\ga_7/\ga_1\eq13$. Since $\ga_1$ is associated with $x^4y^4$ and $13\scd4\eq52$ with $[52/9]\eq5$, $[52/7]\eq7$, we get a division by $9^5\scd7^7\eq48629390607$ during the calculation of $\dd_t^{12}$, but this seems too large for Singular. It may be difficult to calculate this example without replacing $f_1$ with $\tfrac{1}{9}\1x^9+\tfrac{1}{7}\1y^7$ in order to avoid the above division. After a computer calculation using C (where the computation itself takes less than one second), we can conclude after a substitution that the subspace $V(f_1,\tfrac{92}{63})$ on which $\tfrac{92}{63}$ is unshifted is given by
$$\aligned u_3&=\tfrac{44}{3}u_1^4+{\scriptstyle 4}\,u_1u_2,\q
u_4=\tfrac{748}{5}u_1^6+\tfrac{176}{3}u_1^3u_2-{\scriptstyle 16}\1u_1^2u_3+{\scriptstyle 2}\1u_2^2,\\
u_7&\eq-\tfrac{1444507328}{14175}u_1^{13}-\tfrac{20975504}{945}u_1^{10}u_2-\tfrac{11696}{15}u_1^7u_2^2+\tfrac{676}{9}u_1^4u_2^3\\&+\tfrac{6424}{15}u_1^5u_5+\tfrac{8}{3}u_1u_2^4+{\scriptstyle 48}\1u_1^2u_2u_5+\tfrac{176}{3}u_1^3u_6+{\scriptstyle 4}\1u_2u_6.\endaligned$$
This can be obtained by looking at [$*$,7,3] for $*\eq1,2,3$ in the output of the code in \hl{A2}{A.2} of Appendix after typing 9~7~7. Setting $u_j\eq 1$ for $j\eq1,2,5,6$, we obtain that $u_3\eq\tfrac{56}{3}$, $u_4\eq{-}\tfrac{442}{5}$, $u_7\eq{-}\tfrac{1761450728}{14175}$, where the last numerator is quite close to the {\it integer limit\1} in Singular: $2147483647\,({=}\,2^{31}\mi 1)$. By a computation using ``bernstein" in Singular, however, it is rather impressive to see that $\tfrac{92}{63}$ is the only {\it unshifted shiftable\1} root up to sign of the BS polynomial of
$$f\eq\tfrac{1}{9}x^9\pl\tfrac{1}{7}y^7\pl x^4y^4\pl x^3y^5\pl\tfrac{56}{3}x^7y^2-\tfrac{442}{5}x^6y^3\pl x^5y^4\pl x^4y^5\mi\tfrac{1761450728}{14175}x^7y^3.$$
There are six shifted roots up to sign between $\tfrac{92}{63}$ and $\tfrac{79}{63}\eq\alt_f{+}1$, since $\tfrac{92}{63}$ corresponds to $j\eq7$. Hence ${\rm SR}(f,\tfrac{92}{63})\eq\tfrac{1}{2}$ and ${\rm SD}(f,\tfrac{92}{63})\eq\tfrac{2}{9}$. The distribution of roots up to sign is as below (see \hl{2.1}{2.1} for the notation):
$$\setlength{\unitlength}{0.15cm}
\begin{picture}(100,3.5)
\multiput(73,2.5)(7,0){1}{\circle{.7}}
\multiput(75,2.5)(7,0){2}{\circle{.7}}
\multiput(70,2.5)(7,0){4}{\circle{.7}}
\multiput(72,2.5)(7,0){5}{\circle{.7}}
\multiput(6,2.5)(7,0){8}{\circle*{.7}}
\multiput(15,2.5)(7,0){8}{\circle*{.7}}
\multiput(24,2.5)(7,0){7}{\circle*{.7}}
\multiput(33,2.5)(7,0){6}{\circle*{.7}}
\multiput(42,2.5)(7,0){4}{\circle*{.7}}
\multiput(51,2.5)(7,0){3}{\circle*{.7}}
\multiput(9,2.5)(7,0){5}{\circle*{.7}}
\multiput(7,2.5)(7,0){4}{\circle*{.7}}
\put(12,2.5){\circle*{.7}}
\put(10,2.5){\circle*{.7}}
\put(82,2.5){\circle*{.7}}
\put(4.8,0){$\scriptscriptstyle\frac{16}{63}$}
\put(66.8,0){$\scriptscriptstyle\frac{78}{63}$}
\put(81,0){$\scriptscriptstyle\frac{92}{63}$}
\put(69,0){$\scriptscriptstyle\frac{80}{63}$}
\put(98.8,0){$\scriptscriptstyle\frac{110}{63}$}
\end{picture}$$
\par\htt{2.5}{}\msn
{\bf 2.5.~Example V.} Let $f\eq\tfrac{1}{10}\1x^{10}\pl\tfrac{1}{10}\1y^{10}\pl u_1x^3y^8\pl u_6x^8y^3\pl u_9x^6y^6\pl u_{17}x^7y^7$ ($u_i\ins\C)$. By a calculation similar to the weighted homogeneous case, we can see that the unshift condition for the root $\tfrac{16}{10}$ up to sign is given by
$$u_{17}\eq 192(u_1u_6^3\pl u_1^3u_6)\mi 64\1u_1u_6u_9,\q u_9\eq 4\1u_1^2\pl 4\1u_6^2,$$
looking at [1;7,7], [5;7,7] in the output of the code in \hl{A2}{A.2} of Appendix after typing 10~10~18, where we use a variation of Proposition~\hl{2}{2} assuming $u_1u_6\nes 0$. For $u_1\eq u_6\eq 1$, we then obtain that $u_9\eq8$, $u_{17}\eq{-}128$ as is written at the end of the introduction. By a similar argument, we can see, assuming $u_1u_6\nes 0$, that the root $\tfrac{15}{10}$ up to sign is shifted if
$$u_9\nes 4\1u_1^2\pl\tfrac{7}{2}u_6^2,\q u_9\nes\tfrac{7}{2}\1u_1^2\pl 4\1u_6^2,$$
looking at [2;6,7], [3;7,6] in the output of the code in \hl{A2}{A.2} of Appendix. It is easy to see that the roots $\tfrac{13}{10}$, $\tfrac{14}{10}$ up to sign are unshifted. The argument is rather different from the weighted homogeneous case with condition (\hl{M1}{M1}) satisfied.
\par\htt{R2.5}{}\msn
{\bf Remark~2.5.} Setting $f\eq\tfrac{1}{12}\1x^{12}\pl\tfrac{1}{12}\1y^{12}\pl x^4y^{10}\pl x^{10}y^4\pl 10\1x^8y^8$, its BS polynomial has roots $\bl\{\tfrac{2}{12},\dots,\tfrac{18}{12}\br\}\cup\bl\{\tfrac{20}{12}\br\}$ up to sign by Singular. This polynomial has fewer terms than the above example.
\sk
If we set $f\eq\tfrac{1}{14}\1x^{14}\pl\tfrac{1}{14}\1y^{14}\pl x^4y^{11}\pl x^{11}y^4\pl 11\1x^8y^8\mi\tfrac{968}{3}\1x^9y^9\pl 19360\1x^{10}y^{10}$, then its BS polynomial has roots $\bl\{\tfrac{2}{14},\dots,\tfrac{19}{14}\br\}\cup\bl\{\tfrac{22}{14}\br\}$ up to sign, hence $\tfrac{20}{14},\tfrac{21}{14}$ are not, by Singular. These were found by studying pictures like below and modifying the code in \hl{A2}{A.2} of Appendix.
$$\h{$\setlength{\unitlength}{2.5mm}
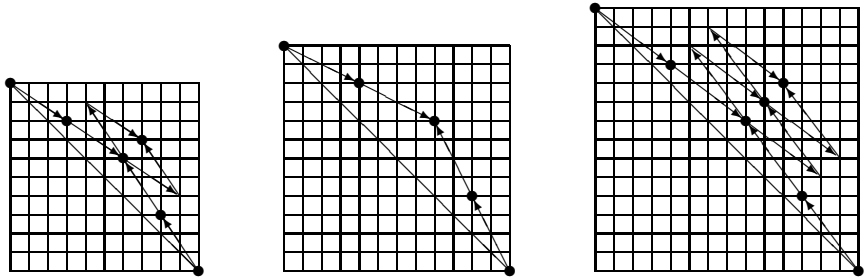
\begin{picture}(14,15)
\multiput(0,0)(1,0){11}{\line(0,1){10}}
\multiput(0,0)(0,1){11}{\line(1,0){10}}
\put(0,10){\line(1,-1){10}}
\put(10,0){\circle*{.5}}
\put(0,10){\circle*{.5}}
\put(8,3){\circle*{.5}}
\put(3,8){\circle*{.5}}
\put(6,6){\circle*{.5}}
\put(7,7){\circle*{.5}}
\multiput(10,0)(-2,3){3}{\vector(-2,3){1.9}}
\multiput(0,10)(3,-2){3}{\vector(3,-2){2.85}}
\put(9,4){\vector(-2,3){1.9}}
\put(4,9){\vector(3,-2){2.85}}
\end{picture}$
$\setlength{\unitlength}{2.5mm}
\begin{picture}(16,15)
\multiput(0,0)(1,0){13}{\line(0,1){12}}
\multiput(0,0)(0,1){13}{\line(1,0){12}}
\put(0,12){\line(1,-1){12}}
\put(12,0){\circle*{.5}}
\put(0,12){\circle*{.5}}
\put(10,4){\circle*{.5}}
\put(4,10){\circle*{.5}}
\put(8,8){\circle*{.5}}
\multiput(12,0)(-2,4){2}{\vector(-1,2){1.9}}
\multiput(0,12)(4,-2){2}{\vector(2,-1){3.8}}
\end{picture}$
$\setlength{\unitlength}{2.5mm}
\begin{picture}(14.2,15)
\multiput(0,0)(1,0){15}{\line(0,1){14}}
\multiput(0,0)(0,1){15}{\line(1,0){14}}
\put(0,14){\line(1,-1){14}}
\put(14,0){\circle*{.5}}
\put(0,14){\circle*{.5}}
\put(11,4){\circle*{.5}}
\put(4,11){\circle*{.5}}
\put(8,8){\circle*{.5}}
\put(9,9){\circle*{.5}}
\put(10,10){\circle*{.5}}
\multiput(14,0)(-3,4){3}{\vector(-3,4){2.85}}
\multiput(0,14)(4,-3){3}{\vector(4,-3){3.8}}
\multiput(12,5)(-3,4){2}{\vector(-3,4){2.85}}
\multiput(5,12)(4,-3){2}{\vector(4,-3){3.8}}
\put(13,6){\vector(-3,4){2.85}}
\put(6,13){\vector(4,-3){3.8}}
\end{picture}$}$$
Here all the arrows are not written and there are much more ``paths" (except in the second). This is closely related to Remark~\hl{R1.2}{1.2}. It seems more difficult to construct an example without symmetry.
\par\htt{2.6}{}\msn
{\bf 2.6.~Example VI.} Set $f\eq\tfrac{1}{6}\1x^6y\pl\tfrac{1}{5}\1xy^5\pl x^5y^2\pl\tfrac{60}{29}\1x^4y^3\pl\tfrac{1320}{841}\1x^3y^4\mi\tfrac{9504000}{594823321}\1x^4y^4$. One can get the coefficients by typing lr6~5~5 for the code in \hl{A2}{A.2} of Appendix and looking at [$*$;4,4] for $*\eq1,2,3$ in the output as is explained before the code. (Here ``mixed fractions" may be used to avoid integer overflow.) The weights of parameters multiplied by $5\1{\cdot}\16{-1}\eq29$ and the exponents of the corresponding monomials are as follows:
$$\scalebox{0.85}{$\begin{array}{cccccc}3&7&11&&\\&2&6\\&&1 \end{array}
\begin{array}{ccccc}3,4&4,4&5,4\\&4,3&5,3\\&&5,2\end{array}$}$$
We have by Singular
$$R_f=\bl\{\tfrac{9}{29},\dots,\tfrac{37}{29}\br\}\cup\bl\{\tfrac{45}{29}\br\}\setminus\bl\{\tfrac{16}{29}\br\}.$$
Hence ${\rm SR}(f,\tfrac{45}{29})=\tfrac{2}{3}$ and ${\rm SD}(f,\tfrac{45}{29})=\tfrac{8}{29}$. The distribution of roots up to sign is as below
$$\hskip-2cm\h{$\setlength{\unitlength}{0.3cm}
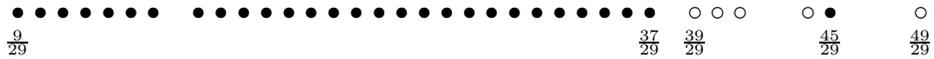
\begin{picture}(50,2.5)
\multiput(9,1.5)(1,0){7}{\circle*{.5}}
\multiput(17,1.5)(1,0){21}{\circle*{.5}}
\multiput(39,1.5)(1,0){3}{\circle{.5}}
\put(45,1.5){\circle*{.5}}
\put(44,1.5){\circle{.5}}
\put(49,1.5){\circle{.5}}
\put(8.4,0){$\scriptscriptstyle\frac{9}{29}$}
\put(36.4,0){$\scriptscriptstyle\frac{37}{29}$}
\put(38.4,0){$\scriptscriptstyle\frac{39}{29}$}
\put(44.4,0){$\scriptscriptstyle\frac{45}{29}$}
\put(48.4,0){$\scriptscriptstyle\frac{49}{29}$}
\end{picture}$}$$
\sk
Setting $\,f\eq\tfrac{1}{7}\1x^7y\pl\tfrac{1}{6}\1xy^6\pl x^6y^2\pl u_2x^5y^3\pl u_3x^4y^4\pl u_4x^3y^5\pl u_6x^5y^4\pl u_7x^4y^5\pl u_9x^5y^5\,$ with $\,u_j\ins\C$ appropriate, we should have ${\rm SR}(f,\tfrac{66}{41})\eq\tfrac{4}{5}$ with $41\eq7\1{\cdot}\1 6{-}1$, although it is difficult to determine $u_6,u_7,u_9$ by the integer overflow problem. The other $u_j$ are given rather easily, and $\tfrac{66}{41}$ is a root up to sign of $b_f(s)$ for {\it some\1} $u_6,u_7,u_9\ins\C$, but it is rather complicated to see that the other shiftable roots are really shifted. These examples are extended to Conjecture~\hl{Cn1}{1}.
\par\htt{R2.6}{}\msn
{\bf Remark~2.6.} Let $f\eq\tfrac{1}{7}\1x^7\pl\tfrac{1}{5}\1xy^5\pl x^6y\pl\tfrac{13}{5}\1x^5y^2\pl\tfrac{221}{75}\1x^4y^3\pl\tfrac{547859}{196875}\1x^5y^3$. This is obtained by typing cr7~5~5 in \hl{A2}{A.2} of Appendix and looking at [$*$;5,3] as is explained before the code. By Singular we see that $\tfrac{54}{35}$ is the unique unshifted shiftable root up to sign of $b_f(s)$ with ${\rm SR}(f,\tfrac{54}{35})=\tfrac{2}{3}$ and ${\rm SD}(f,\tfrac{54}{35})=\tfrac{2}{7}$. The distribution of roots up to sign is as below (see \hl{2.1}{2.1} for the notation):
$$\h{$\setlength{\unitlength}{0.3cm}
\begin{picture}(50,2.5)
\multiput(1,1.5)(1,0){4}{\circle*{.5}}
\multiput(6,1.5)(1,0){3}{\circle*{.5}}
\multiput(11,1.5)(1,0){4}{\circle*{.5}}
\multiput(16,1.5)(1,0){4}{\circle*{.5}}
\multiput(21,1.5)(1,0){9}{\circle*{.5}}
\multiput(31,1.5)(1,0){4}{\circle*{.5}}
\multiput(44,1.5)(1,0){1}{\circle*{.5}}
\multiput(37,1.5)(1,0){3}{\circle{.5}}
\put(43,1.5){\circle{.5}}
\put(49,1.5){\circle{.5}}
\put(.4,0){$\scriptscriptstyle\frac{11}{35}$}
\put(33.4,0){$\scriptscriptstyle\frac{44}{35}$}
\put(36.4,0){$\scriptscriptstyle\frac{47}{35}$}
\put(43.4,0){$\scriptscriptstyle\frac{54}{35}$}
\put(48.4,0){$\scriptscriptstyle\frac{59}{35}$}
\end{picture}$}$$
\par\htt{2.7}{}\msn
{\bf 2.7.~Example VII.} Let $f_1\eq x^7{+}\1y^5{+}z^3$. After $x^5{+}\1y^4{+}z^3$ and $x^7{+}\1y^4{+}z^3$, this is the third simplest example in the three variable case. The weights of parameters multiplied by 105 and the exponents of the corresponding monomials are as follows:
$$\scalebox{0.8}{$\begin{array}{cccccccccc}&&&&5\\
&&&11&26\\
&2&17&32&47\\
8&23&38&53&68&12&&&\\
&&3&18&33\end{array}
\begin{array}{cccccccccc}&&&&5,0,1\\
&&&4,1,1&5,1,1\\
&2,2,1&3,2,1&4,2,1&5,2,1\\
1,3,1&2,3,1&3,3,1&4,3,1&5,3,1&5,2,0\\
&&3,3,0&4,3,0&5,3,0\end{array}$}$$
We study the shift of the root $\tfrac{202}{105}$ up to sign corresponding to $j\eq10$ with $\ga_{10}\eq\tfrac{26}{105}$, where $\alt_f{+}1\eq\tfrac{176}{105}$. By a computer calculation using C, we conclude after a substitution that the subspace $V(f_1,\tfrac{202}{105})$ on which $\tfrac{202}{105}$ is unshifted is defined by
$$\aligned u_3&\eq{\scriptstyle 2}\1u_1u_2,\q u_5\eq{-}{\scriptstyle 7}\1u_1^4u_2\pl{\scriptstyle 4}\1u_1^3u_3\pl{\scriptstyle 14}\1u_1u_2^3\mi{\scriptstyle 6}\1u_2^2u_3\pl{\scriptstyle 2}\1u_2u_4\mi u_5,\\
u_{10}&\eq{-}\tfrac{374}{27}u_1^{13}\pl{\scriptstyle 162}\1u_1^{10}u_2^2\mi{\scriptstyle 678}\1u_1^7u_2^4\pl\tfrac{4176}{5}u_1^4u_2^6\pl\tfrac{442}{9}u_1^9u_4\mi\tfrac{2223}{5}u_1u_2^8\mi{\scriptstyle 285}\1u_1^6u_2^2u_4\\
&\q\pl{\scriptstyle 366}\1u_1^3u_2^4u_4\pl{\scriptstyle 26}\1u_1^7u_6\mi{\scriptstyle 91}\1u_2^6u_4\mi{\scriptstyle 42}\1u_1^5u_4^2\mi{\scriptstyle 142}\1u_1^4u_2^2u_6\pl{\scriptstyle 60}\1u_1^2u_2^2u_4^2\pl{\scriptstyle 158}\1u_1u_2^4u_6\\
&\q\mi{\scriptstyle 28}\1u_1^3u_4u_6\mi{\scriptstyle 4}\1u_1^3u_2u_7\mi{\scriptstyle 7}\1u_1^4u_8\pl\tfrac{14}{3}u_1u_4^3\pl{\scriptstyle 18}\1u_2^2u_4u_6\pl{\scriptstyle 14}\1u_2^3u_7\pl{\scriptstyle 18}\1u_1u_2^2u_8\\
&\q\mi{\scriptstyle 6}\1u_1u_6^2\pl{\scriptstyle 2}\1u_4u_8\pl{\scriptstyle 2}\1u_2u_9.\endaligned$$
It does not necessarily seem easy to determine the BS polynomial of $f$ if we put $u_i\eq1$ for any $j\ins\{1,\dots,9\}\stm\{3,5\}$. Setting $u_1\eq u_2\eq 1$ and $u_j\eq0$ for $j\eq4,6,7,8,9$, we get
$$f\eq\tfrac{1}{7}x^7\pl\tfrac{1}{5}y^5\pl\tfrac{1}{3}z^3\pl x^2y^2z\pl x^3y^3\pl 2x^5z\pl 3x^4yz\mi\tfrac{18799}{135}x^5yz.$$
Its BS polynomial has the unique unshifted shiftable root $\tfrac{202}{105}$ up to sign according to Singular. We have ${\rm SR}(f,\tfrac{202}{105})\eq\tfrac{9}{16}$, ${\rm SD}(f,\tfrac{202}{105})\eq\tfrac{29}{105}$, and the distribution of roots up to sign is as below (see \hl{2.1}{2.1} for the notation):
$$\setlength{\unitlength}{0.9mm}
\begin{picture}(180,8)
\multiput(1,4.5)(15,0){6}{\circle*{0.9}}
\multiput(22,4.5)(15,0){6}{\circle*{0.9}}
\multiput(43,4.5)(15,0){5}{\circle*{0.9}}
\multiput(64,4.5)(15,0){3}{\circle*{0.9}}
\multiput(36,4.5)(15,0){5}{\circle*{0.9}}
\multiput(57,4.5)(15,0){4}{\circle*{0.9}}
\multiput(78,4.5)(15,0){2}{\circle*{0.9}}
\multiput(99,4.5)(15,0){1}{\circle*{0.9}}
\multiput(118,4.5)(15,0){1}{\circle{0.9}}
\multiput(109,4.5)(15,0){3}{\circle{0.9}}
\multiput(111,4.5)(15,0){1}{\circle{0.9}}
\multiput(117,4.5)(15,0){1}{\circle{0.9}}
\multiput(108,4.5)(15,0){4}{\circle{0.9}}
\multiput(114,4.5)(15,0){5}{\circle{0.9}}
\multiput(13,4.5)(15,0){1}{\circle*{0.9}}
\multiput(4,4.5)(15,0){3}{\circle*{0.9}}
\multiput(6,4.5)(15,0){1}{\circle*{0.9}}
\multiput(12,4.5)(15,0){1}{\circle*{0.9}}
\multiput(3,4.5)(15,0){4}{\circle*{0.9}}
\multiput(9,4.5)(15,0){5}{\circle*{0.9}}
\multiput(132,4.5)(15,0){1}{\circle*{0.9}}
\put(100.25,0){$\scriptscriptstyle\frac{173}{105}$}
\put(105.25,0){$\scriptscriptstyle\frac{178}{105}$}
\put(129.3,0){$\scriptscriptstyle\frac{202}{105}$}
\put(171.3,0){$\scriptscriptstyle\frac{244}{105}$}
\end{picture}$$
\par\htt{2.8}{}\msn
{\bf 2.8.~Example VIII.} Let $f_1\eq x^4y\pl y^4z\pl xz^3$. This is a polynomial of nearly BP loop type, where $\mu_{f_1}\eq48$ with $T^{49}\eq1$, see Remark~\hl{R1.1e}{1.1e}. The Jacobian ring is spanned by $x^iy^jz^k$ ($i,j\slt4,k\slt3$) with
weighted degree function given by
$$\ell(i,j,k)\eq(10i{+}9j{+}13k)/49.$$
The modality of $f_1$ (that is, $|J|$) is 17. The minimal, maximal, and the second largest spectral numbers are respectively $\tfrac{32}{49}$, $\tfrac{115}{49}$, $\tfrac{106}{49}$. We see that the last number subtracted by the minimal exponent $\tfrac{32}{49}$ and $\ga_1\eq\tfrac{1}{49}$ added by 1 are associated respectively with the monomials $x^3y^2z^2$ and $xy^3z$, since $\ell(3,2,2)\eq\tfrac{74}{49}$, $\ell(1,3,1)\eq\tfrac{50}{49}$. Observe that $x^3y^2z^2$ is {\it not over\1} $xy^3z$. Setting
$$\aligned f&\eq x^4y\pl y^4z\pl xz^3\pl xy^3z\pl u_2x^2y^2z\pl u_3x^3yz\pl u_4xy^2z^2\pl u_5x^2yz^2\\ &\q\pl u_6x^3z^2\pl u_7x^3y^2z\pl u_8x^2y^2z^2\pl u_9x^3yz^2\pl u_{10}x^3y^2z^2,\endaligned$$
with $u_i\in\C^*$ appropriate (where monomials which are not under $x^3y^2z^2$ are omitted), it is expected that $\tfrac{106}{49}$ is the unique unshifted shiftable root up to sign of $f$ with ${\rm SR}(f,\tfrac{106}{49})\eq\tfrac{15}{17}$. Here the denominators of $u_8$, $u_9$ may be close to $49^{15}$ and $49^{16}$, calculating $(\ell(2,2,2){-}1)/\ga_1$, etc. It is then quite nontrivial to show that the other shiftable roots are really shifted as in Conjecture~\hl{Cn1}{1}. Note that all the shiftable roots are shifted if all the $u_i$ vanish, see Corollary~\hl{C4}{4}. This can be verified by using Singular for this example.
\par\htt{R2.8a}{}\msn
{\bf Remark~2.8a.} It is not easy to compute the $g_{k,l}^{(1)}$ for polynomials of nearly BP loop type for $n\eq 3$, since the action of $\dd_t$ is always associated with the division by $abc{+}1$ which is very often a (square of) relatively large prime. This cannot be avoided by adding appropriate coefficients to $f_1$ as in the BP type case. The denominators tend to be huge even in the case $(a,b,c)\eq(4,3,3)$ or $(5,3,2)$. There is, however, a calculable nontrivial example
\htt{2.8.1}{}
$$f\eq x^6y\pl y^2z\pl xz^2\pl x^5z\mi\tfrac{17}{25}\1x^3yz\pl\tfrac{305966}{1953125}\1x^4yz.
\leqno(2.8.1)$$
Here $\mu_f\eq24$, $R_f\eq\{\tfrac{21}{25},\dots,\tfrac{45}{25}\}\cup\{\tfrac{51}{25}\}\stm\{\tfrac{26}{25}\}$ with ${\rm SR}(f,\tfrac{51}{25})\eq\tfrac{1}{2}$, ${\rm SD}(f,\tfrac{51}{25})\eq\tfrac{6}{25}$ by Singular (and $\Rt_f\eq R_f\stm\{1\}$ in this case).
\par\htt{R2.8b}{}\msn
{\bf Remark~2.8b.} If the last two terms of $f$ in (\hl{2.8.1}{2.8.1}) are omitted, we get that
$$R_f\eq\{\tfrac{21}{25},\dots,\tfrac{45}{25}\},$$ see Corollary~\hl{C4}{4}. This gives an example such that the roots of BS polynomial are {\it consecutive\1} with common denominator the order of the monodromy and moreover every monodromy eigenvalue of $f$ has {\it multiplicity\1} 1, see also Remark~\hl{R1.1e}{1.1e}.
\bs\bs
\vbox{\centerline{\bf Appendix: Sample codes}
\bsn
In this Appendix we give some sample codes to compute the $g_{k,l}^{(1)}$ and the bistable subsets.}
\par\htt{A1}{}\msn
{\bf A.1.~Sample code for Singular.} For the convenience of the reader we note here a sample code to calculate Example~\hl{2.3}{III} using Singular. Since it is written in a condensed way, it may be better to add line breaks appropriately after copying and pasting it in a text file. One may modify a,b as long as a+b\,$\less$\,15 and (a,b)=1. (If ``Division Error" appears, one has to increase the size of the vector iv.) This code cannot be applied to $f_1\eq x^9{+}y^7$, although it works at least for (a,b)\,=\,(7,6), (9,4), (7,5). (Please verify whether the list of weights in the last line is correct.)
\ms
\vbox{\tiny\sf\pv@ring R = 0, (u_1,u_2,u_3,u_4,u_5,u_6,u_7), ds;@
\pv@int a,b,rs,num,i,j,wd,n,p,iq,ir,jq,jr,kk,e,li,lj,rs,maxp,MA,MB,NuM;a=8;b=7;@
\pv@int k,q,od,di,ip,jp,wp,ie,je,we,pp,rp,MM,mxn,mxk,maxdiv,am,bm,m,mm;poly Sub;m=a*b;@
\pv@vector iv=[1,1/2,1/3,1/4,1/5,1/6,1/7,1/8,1/9,1/10,1/11,1/12];intvec wt=0,0,0,0,0,0,0;@
\pv@am=a-1;bm=b-1;mm=am*bm;if(a+b>14){rs=5;}else{rs=4;}mxn=7;MA=2000;MB=am*bm;NuM=1000;@
\pv@intmat Co[mxn][2];intmat Nu[NuM][NuM];matrix M[MA][MB];intmat O[MA][MB];@
\pv@vector va=[u_1,u_2,u_3,u_4,u_5,u_6,u_7];mxk=2*(m-a-b);p=1;for(k=m+1;k<mxk&&p<=mxn;k++)@
\pv@{for(i=1;i<a-1&&k>b*i;i++){j=(k-b*i)div a;if((k-b*i)%a==0&&j<b-1){Co[p,1]=i;Co[p,2]=j;@
\pv@wt[p]=k-m;p++;}}}num=p-1;rp=wt[num]div wt[1];MM=wt[num]+m*rp;p=0;maxdiv=0;for(wd=MM;@
\pv@wd>=0;wd--){n=wd div b;for(i=0;i<=n;i++){j=(wd-b*i)div a;if((wd-b*i)%a==0){p++;@
\pv@Nu[i+1,j+1]=p;ir=i%a;iq=i div a;jr=j%b;jq=j div b;li=1;for(e=1;e<=iq;e++){li=li*(i-e*@
\pv@a+1);}lj=1;for(e=1;e<=jq;e++){lj=lj*(j-e*b+1);}if(ir!=a-1&&jr!=b-1){M[p,ir+am*jr+1]=@
\pv@li*lj*iv[a]^iq*iv[b]^jq;O[p,ir+am*jr+1]=iq+jq;}for(q=1;q<=num;q++){ip=i+Co[q,1];jp=j+@
\pv@Co[q,2];wp=b*ip+a*jp;if(wp<=MM){pp=Nu[ip+1,jp+1];for(e=1;e<=mm;e++){ie=(e-1)%am;je=@
\pv@(e-1)div am;we=b*ie+a*je;if(M[pp,e]!=0&&ie+je<a+b-rs&&we>m){od=O[pp,e];di=we+m*(od-1)@
\pv@-wd;if(di<=wt[num]){Sub=M[pp,e]*wt[q]*iv[di]*va[q];if(O[p,e]==od-1){if(di>maxdiv)@
\pv@{maxdiv=di;}M[p,e]=M[p,e]-Sub;}if(O[p,e]>od-1||(O[p,e]<od-1&&M[p,e]==0)){if(di>maxdiv)@
\pv@{maxdiv=di;}M[p,e]=0-Sub;O[p,e]=od-1;}}}}}}}}} maxp=p;if(size(iv)<maxdiv){sprintf(@
\pv@"Division Error %s",maxdiv);} for(i=1;i<=maxp;i++){for(e=1;e<=mm;e++){ie=(e-1)%am;je@
\pv@=(e-1)div am;we=b*ie+a*je;if(O[i,e]==-1&&ie+je<a+b-rs&&we>m){sprintf("[%s;%s,%s]:",@
\pv@maxp+1-i,(e-1)%am,(e-1)div am);M[i,e];}}} sprintf("wt=%s",wt);@}
\msn
In case one needs a computation for $\tfrac{1}{a}\1x^a{+}\tfrac{1}{b}\1y^b$ (instead of $x^a{+}y^b$) as in Example~\hl{2.4}{IV}, one can do it by removing {\small\verb@*iv[a]^iq*iv[b]^jq@} i.
\par\htt{A2}{}\msn
{\bf A.2.~Sample code using C.} There is also a sample code by C as below. This is possible, since the algorithm is extremely simple. One may copy and past it in a text file and compile it using gcc in Unix or clang in Mac, etc. In the first two lines, {\it line breaks\1} must be inserted before \# and ``long long". Do not forget to remove page numbers, etc. (It may be necessary to replace ' with a {\it character from keyboard\1} if ``Preview" in Mac is used.) When one runs ./a.out, one is asked to type three numbers (and press the return key), which are 9~7~7 and 10~10~18 in the case of Examples~\hl{2.4}{IV} and \hl{2.5}{V}. The last number is the number of variables of the parameter space in which one is interested. (It is not assumed that the exponents are mutually prime. Note however that the code can calculate only the {\it strictly negative grading part\1} of the parameter space of the miniversal $\mu$-constant deformation unless the exponents are mutually prime.) One can type \,n\, before these numbers if one needs a computation for $f_1\eq x^a{+}y^b$ instead of $f_1\eq\tfrac{1}{a}\1x^a{+}\tfrac{1}{b}\1y^b$. So one may enter n7~5~3, n9~4~3, n7~6~5, n8~7~7 for examples in \hl{2.1}{2.1}--\hl{2.3}{3}. 
\ms
\vbox{\fontsize{7pt}{2.5mm}\sf\pv@#include<stdio.h> #include<stdlib.h> #define MN 3500 #define MR 85 #define PN 21 #define NuM 5000@
\pv@long long L0,L1,L2,L3;int P1[PN],P2[PN],P3[PN],PP[PN],AA[PN],BB[PN],DI[PN],Dia[PN][PN],Cn[PN][PN],@
\pv@J[PN][3],JW[PN][PN],Ad[MN][2],num,maxn,maxp,Md,e,p,cc,M[MN][PN][MR][2*PN],Ord[MN][PN],Mpr[2*PN],@
\pv@N[MN][PN],Nu[NuM][NuM],Prm[PN],PRMS[10*PN],Cw[PN][PN],O[2*PN][2*PN],S[2*PN][MR][2*PN],Nn[2*PN],@
\pv@R[2*PN],SG[MR],W[MR],WW[MR],a,b,m,amm,bmm,ss,nump,bg,i,j,ir,jr,gcd,add,bdd,idd,jdd,iq,jq,ee,n;@
\pv@long long li,lj,L[MN][PN][MR];int GCD(int x,int y);int gtd(void);void prf(long long LLL);@
\pv@long long summ(void);long long pfac(long long x);void SS(void);void AAA(void);void BBB(void);@
\pv@void prn(int r);void s(int r);int main(void) {int c,k,l,nn,nnn,q,wd,rp,fl,mw,am,bm,MM,ad,bd,kk,@
\pv@od,di,aa,bb,ii,jj,ip,jp,wp,ie,je,we,pp;long long ll,lll;for(i=2;i<10*PN;i++){PRMS[i]=1;}for(i=2;@
\pv@i*i<10*PN;i++){for(j=2;i*j<10*PN;j++){PRMS[i*j]=0;}}for(i=j=0;i<10*PN&&j<PN;i++){if(PRMS[i]==1)@
\pv@{Prm[j]=i;j++;}}printf("Exponents and Number of members of J: ");for(cc=getchar();cc<'0'||cc>'9';@
\pv@cc=getchar()){if(cc=='c'){Md=1;}if(cc=='l'){Md=2;}if(cc=='n'){n=1;}if(cc=='r'){bg=1;}if(cc=='s')@
\pv@{ss=1;}if(cc=='w'){ss=2;}}if(n==1){Md=Md+10;n=0;}a=aa=gtd();b=bb=gtd();if(ss!=1)num=gtd();am=a-1;@
\pv@bm=b-1;if(a>PN||b>PN){printf("Too big exponents\n"),exit(1);}if(Md%10==1){aa=am;}else if(Md%10==@
\pv@2){aa=am;bb=bm;}pfac(a);for(i=0;i<PN;i++){AA[i]=PP[i];}pfac(b);for(i=0;i<PN;i++){BB[i]=PP[i];}@
\pv@gcd=GCD(aa,bb);if(Md%10==0){ad=a/gcd;amm=am;bmm=bm;bd=b/gcd;m=a*bd;}else if(Md%10==1){ad=am/gcd;@
\pv@amm=a;bmm=bm;bd=b/gcd;m=a*bd;}else if(Md%10==2){ad=am/gcd;amm=a;bmm=b;bd=bm/gcd;m=(a*b-1)/gcd;}@
\pv@for(i=0;i<PN;i++){Dia[i][i]=1;}nump=PN;for(p=0,k=m+1;k<2*m&&p<nump;k++){for(i=1;i<amm&&k>bd*i;@
\pv@i++){if((k-bd*i)%ad==0&&(j=(k-bd*i)/ad)<bmm){J[p][0]=k-m;J[p][1]=i;J[p][2]=j;W[k-m]=Cn[i][j]=@
\pv@p+1;Cw[i][j]=k-m;lll=pfac(J[p][0]);if(lll!=1){printf("Wt Error %d\n",J[p][0]),exit(1);}for(l=0;@
\pv@l<PN;l++){JW[p][l]=PP[l];}p++;}}}if(ss>0){nump=p;SS();if(ss==1)exit(1);}if(num>p||num>PN){printf(@
\pv@"num=%d is replaced by %d!\n",num,p);num=p;}if(bg>0){bg=num-1;}rp=J[num-1][0]/J[0][0];MM=J[num-1][0]@
\pv@+m*rp;printf("Weights and exponents:\n");for(i=0;i<num;i++){printf("%d (%d,%d), ",J[i][0],J[i][1]@
\pv@,J[i][2]);}printf("\nf=");if(Md<3){printf("(1/%d)*",a);}printf("x^%d",a);if(Md%10==2){printf(@
\pv@"*y");}printf("+");if(Md<3){printf("(1/%d)*",b);}if(Md%10==1||Md%10==2){printf("x*");}printf("y^%d",@
\pv@b);for(p=0;p<num;p++){printf("+u_%d*x^%d*y^%d",p+1,J[p][1],J[p][2]);}Nu[0][0]=Nu[1][0]=Nu[0][1]=@
\pv@Nu[2][0]=Nu[1][1]=Nu[0][2]=-1;for(p=wd=0;wd<=MM;wd++){nnn=wd/bd;for(i=0;i<=nnn;i++){j=(wd-bd*i)/ad;@
\pv@if((wd-bd*i)%ad==0&&Nu[i][j]<0){if(p>=MN){printf("Too many monomials: %d\n",p),exit(1);}Nu[i][j]=p+@
\pv@1;Ad[p][0]=i;Ad[p][1]=j;for(k=0;k<num;k++){if(Nu[i+J[k][1]][j+J[k][2]]<=0&&i+J[k][1]<NuM&&j+J[k][2]@
\pv@<NuM){Nu[i+J[k][1]][j+J[k][2]]=-1;}}p++;}}}printf(";\nCalculating %d monomials",p);for(p--;p>=0;p--)@
\pv@{i=Ad[p][0];j=Ad[p][1];wd=ad*j+bd*i;ir=i%a;iq=i/a;jr=j%b;jq=j/b;for(li=e=1;e<=iq;e++){li=li*(i-e*a+@
\pv@1);}for(lj=e=1;e<=jq;e++){lj=lj*(j-e*b+1);}if(Md%10==1){ir=(i-jq)%a;iq=(i-jq)/a;}else if(Md%10==2){@
\pv@iq=(b*i-j)/(a*b-1);jq=(a*j-i)/(a*b-1);ir=i-a*iq-jq;jr=j-iq-b*jq;}ee=Cn[ir][jr]-1;n=N[p][ee];if(ir!=@
\pv@amm&&jr!=bmm&&ee>=0){if(Md%10==0){L[p][ee][n]=pfac(li*lj);for(l=0;l<PN;l++){M[p][ee][n][l]=0-PP[l];@
\pv@}}else if(Md%10==1){AAA();}else if(Md%10==2){BBB();}if(Md%20==10){for(l=0;l<PN;l++){M[p][ee][n][l]@
\pv@=M[p][ee][n][l]+AA[l]*iq+BB[l]*jq;}}Ord[p][ee]=iq+jq;N[p][ee]++;if(N[p][ee]>maxn)maxn=N[p][ee];if(@
\pv@N[p][e]>=MR){printf("Too many terms %d\n",N[p][e]),exit(1);}}for(q=0;q<num;q++){ip=i+J[q][1];jp=j@
\pv@+J[q][2];wp=bd*ip+ad*jp;if(wp<=MM){pp=Nu[ip][jp]-1;for(e=0;e<num;e++){nn=N[pp][e];ie=J[e][1];je=@
\pv@J[e][2];we=bd*ie+ad*je;od=Ord[pp][e];di=we+m*(od-1)-wd;if(nn!=0&&di<=J[num-1][0]){ll=pfac(di);if@
\pv@(ll!=1){printf("Too large di=%d,ll=%lld,we=%d,od=%d,wd=%d\n",di,ll,we,od,wd);}for(l=0;l<PN;l++){@
\pv@DI[l]=PP[l];}if(Ord[p][e]==od-1&&N[p][e]!=0){for(cc=0;cc<nn;cc++){n=N[p][e];for(fl=c=0;c<n&&fl==0;@
\pv@c++){for(l=0;l<num&&M[p][e][c][PN+l]==M[pp][e][cc][PN+l]+Dia[q][l];l++){;}if(l==num){for(l=0;l<PN;@
\pv@l++){P1[l]=M[p][e][c][l];P2[l]=M[pp][e][cc][l]-JW[q][l]+DI[l];}L1=L0=L[p][e][c];L2=0-L[pp][e][cc];@
\pv@L[p][e][c]=summ();for(l=0;l<PN;l++){M[p][e][c][l]=P3[l];}fl++;}}if(fl==0){L[p][e][n]=0-L[pp][e][cc];@
\pv@for(l=0;l<PN;l++){M[p][e][n][l]=M[pp][e][cc][l]-JW[q][l]+DI[l];}for(l=0;l<num;l++){M[p][e][n][l+PN]@
\pv@=M[pp][e][cc][l+PN]+Dia[q][l];}N[p][e]++;if(N[p][e]>maxn)maxn=N[p][e];if(N[p][e]>=MR){printf(@
\pv@"Too many terms %d\n",N[p][e]),exit(1);}}}}if(Ord[p][e]>od-1||N[p][e]==0){for(cc=0;cc<N[pp][e];cc++)@
\pv@{L[p][e][cc]=0-L[pp][e][cc];for(l=0;l<PN;l++){M[p][e][cc][l]=M[pp][e][cc][l]-JW[q][l]+DI[l];}for(l=0;@
\pv@l<num;l++){M[p][e][cc][l+PN]=M[pp][e][cc][l+PN]+Dia[q][l];}}N[p][e]=N[pp][e];Ord[p][e]=od-1;}}}}}}for@
\pv@(i=5;i>=0;i--){for(e=bg;e<num;e++){if(Ord[i][e]==-1){printf(";\n[%d;%d,%d]:",i+1,J[e][1],J[e][2]);for@
\pv@(j=N[i][e]-1;j>=0;j--){for(l=0;l<PN+num;l++){Mpr[l]=M[i][e][j][l];}prf(L[i][e][j]);}}}}printf(";\n");@
\pv@}void prf(long long LLL){int j,k,l,fl;long long De;fl=0;printf("\n");for(De=1,l=0;l<PN;l++){if(Mpr[l]@
\pv@>0){for(j=0;j<Mpr[l];j++){De=De*Prm[l];}}}for(l=0;l<PN;l++){if(Mpr[l]<0){for(j=0;j<0-Mpr[l];j++){LLL@
\pv@=LLL*Prm[l];}}}if(LLL>0&&De>1){printf("+%lld/%lld",LLL,De);fl=1;}else if(LLL>1&&De==1){printf("+%lld"@
\pv@,LLL);fl=1;}else if(LLL<0&&De>1){printf("%lld/%lld",LLL,De);fl=1;}else if(LLL<-1&&De==1){printf("%lld"@
\pv@,LLL);fl=1;}else if(LLL==-1&&De==1){printf("-");}else if(LLL==1&&De==1){printf("+");}else if(LLL==0)@
\pv@{printf("0!!! ");}for(l=0;l<num;l++){if(fl==0&&Mpr[l+PN]!=0){fl=1;}else if(fl!=0&&Mpr[l+PN]!=0){@
\pv@printf("*");}if(Mpr[l+PN]>1){printf("u_%d^%d",l+1,Mpr[l+PN]);}else if(Mpr[l+PN]==1){printf("u_%d",@
\pv@l+1);}}}long long summ(void){int l,r,k;for(l=0;l<PN;l++){if(P1[l]>P2[l]){for(k=1,r=0;r<P1[l]-P2[l];@
\pv@r++){k=k*Prm[l];}L2=L2*k;P3[l]=P1[l];}else if(P1[l]<P2[l]){for(k=1,r=0;r<P2[l]-P1[l];r++){k=k*Prm[l]@
\pv@;}L1=L1*k;P3[l]=P2[l];}else{P3[l]=P2[l];}}L3=L1+L2;if((L1>0&&L2>0&&L3<0)||(L1<0&&L2<0&&L3>0)){printf@
\pv@("\nIntOvFl at p=%d (%d,%d)",p,J[e][1],J[e][2]);}L3=pfac(L3);for(l=0;l<PN;l++){P3[l]=P3[l]-PP[l];}@
\pv@return(L3);}long long pfac(long long x){lldiv_t dvn;int l,flg,fl2,sgn=1;if(x<0){x=-x;sgn=-1;}for(l@
\pv@=0;l<PN;l++){PP[l]=0;}for(flg=0;x>1&&flg==0;){for(l=fl2=0;l<PN&&fl2==0;l++){dvn=lldiv(x,Prm[l]);if@
\pv@(dvn.rem==0){fl2++;if(l>maxp){maxp=l;}}}if(fl2==0){flg++;}else{x=dvn.quot;PP[l-1]++;}}if(sgn<0){x=@
\pv@-x;}return(x);}int gtd(void){int l;while(cc<'0'||'9'<cc){cc=getchar();}for(l=0;'0'<=cc&&cc<='9';cc@
\pv@=getchar()){l=l*10+cc-'0';}return(l);}int GCD(int x,int y){while(x!=0&&y!=0){if(x>y){x=x%y;}else{y@
\pv@=y%x;}}if(x==0){return(y);}else{return(x);}}void SS(void){printf("Replace SS!\n");exit(1);}void@
\pv@AAA(void){printf("Replace AAA!\n");exit(1);}void BBB(void){printf("Replace BBB!\n");exit(1);}@}
\msn
If one replaces the definitions of ``void AAA(void)" and ``void BBB(void)" at the end of the above code with the code below, one can type c or l before the numbers in the case $f_1$ is a polynomial of nearly BP {\it chain\1} or {\it loop\1} type, see \hl{2.6}{2.6} and Remark~\hl{R1.1c}{1.1c}. For a calculation using Singular as in Remark~\hl{R2.3}{2.3}, one can replace the {\it first\1} term ``{\small\sf\verb@-u_k@}" (or ``{\small\sf\verb@-u_k+@}") of the polynomials at ``[$*$,i,j]:" in the output by ``{\small\sf\verb@poly u_k=@}" for a fixed (i,j), and use the expression of $f$ in the output in order to apply ``bernstein(f);" setting ``{\small\sf\verb@poly u_k=1@}" (usually at least for k=1) or ``{\small\sf\verb@0@}" appropriately for the {\small\sf\verb@k@} not appearing as the first term of polynomials in [$*$,i,j], see also \hl{2.6}{2.6}. The value of {\small\sf\verb@u_k@} can be obtained by typing ``{\small\sf\verb@u_k;@}". If one prefers to see only the data for the last member, one can type r.
\ms
\vbox{\fontsize{7pt}{2.5mm}\sf\pv@void AAA(void){int k,l,ii,jj;long ll;if(pfac(jr+1)!=1){printf("pfac error: jr+1=%d\n",jr+1),exit(1);@
\pv@}for(l=0;l<PN;l++){if(PP[l]>BB[l]){P3[l]=0;}else{P3[l]=BB[l]-PP[l];}}gcd=GCD(b,jr+1);bdd=b/gcd;jdd=@
\pv@(jr+1)/gcd;for(l=1,ll=1;l<=iq;l++){ll=ll*(bdd*(i-jq-a*l+1)-jdd);}for(l=0;l<PN;l++){P3[l]=P3[l]*iq;}@
\pv@L[p][ee][n]=pfac(ll*lj);for(l=0;l<PN;l++){M[p][ee][n][l]=P3[l]-PP[l];}if(Md%20==11){for(l=0;l<PN;@
\pv@l++){M[p][ee][n][l]=M[p][ee][n][l]+iq*AA[l];}}}void BBB(void){int k,l,ii,jj;long ll;ii=i;jj=j;for@
\pv@(k=1,ll=1;k<=iq;k++){pfac(jj);for(l=0;l<PN;l++){if(PP[l]>BB[l]){P3[l]=0;}else{P3[l]=BB[l]-PP[l];}}@
\pv@gcd=GCD(b,jj);bdd=b/gcd;jdd=jj/gcd;ll=ll*pfac(bdd*(ii-a+1)-jdd);for(l=0;l<PN;l++){M[p][ee][n][l]=@
\pv@M[p][ee][n][l]+P3[l]-PP[l];}ii=ii-a;jj=jj-1;}for(k=1;k<=jq;k++){pfac(ii);for(l=0;l<PN;l++){if(PP[l]@
\pv@>AA[l]){P3[l]=0;}else{P3[l]=AA[l]-PP[l];}}gcd=GCD(a,ii);add=a/gcd;idd=ii/gcd;ll=ll*pfac(add*(jj-b+@
\pv@1)-idd);for(l=0;l<PN;l++){M[p][ee][n][l]=M[p][ee][n][l]+P3[l]-PP[l];}ii=ii-1;jj=jj-b;}pfac(a*b-1);@
\pv@for(l=0;l<PN;l++){if(Md%20==2){M[p][ee][n][l]=M[p][ee][n][l]+(iq+jq)*(PP[l]-AA[l]-BB[l]);}else{@
\pv@M[p][ee][n][l]=M[p][ee][n][l]+(iq+jq)*PP[l]-jq*AA[l]-iq*BB[l];}}L[p][ee][n]=ll;}@}
\par\htt{RA2}{}\msn
{\bf Remark~A.2.} This algorithm was at first suspected to be too simple to be true. However, many explicit calculations using the codes in \hl{A1}{A.1}--\hl{A2}{2} combined with ``bernstein" in Singular seem to show that it is correct. Please inform us in case one gets a different result using the codes in \hl{A1}{A.1} and \hl{A2}{A.2}. (Recall that one can modify the code in \hl{A1}{A.1} as is noted at the end of it in order to compute examples with coefficients of $f_1$ modified.) Since Gr\"obner bases are not used in the second code, the orderings of monomials are different in general, although they may coincide in simple cases as a consequence of the algorithm. The second code investigates only the necessary monomials, and this is the main reason for which it is very fast. To calculate more complicated examples, one may like to change some parameters defined at the beginning, and this may be possible in case ``virtual memory" is available although the computation might become slow. Anyway it does not seem easy to handle more than 64 bit integers. Since localization by small prime numbers is used in the code, integer overflow may occur at the last stage where the conversion to usual rational numbers is made. These codes are still experimental, and might contain a bug that does not appear for relatively simple examples.
\par\htt{A3}{}\msn
{\bf A.3.~Sample code calculating the bistable subsets.} We note here a sample code to compute the bistable subsets in the two variable case with {\it mutually prime\1} modified exponents, where the {\it modified exponents\1} are $(a{-}1,b)$ and $(a{-}1,b{-}1)$ respectively in the chain and loop type cases. It is designed to be used by replacing the definition of ``void SS(void)" at the end of the code in \hl{A2}{A.2} by it, and typing s before the numbers, for instance s8~7. One can type w (together with the third number) instead of s in case both outputs are needed. (Here the modified exponents must be mutually prime.) It is {\it not\1} necessary to type the letters n, r, c, l, s, w in alphabetical order.
\ms
\vbox{\fontsize{7pt}{2.5mm}\sf\pv@void SS(void){int c,i,j,e,k,l,n,p,q,fl,mw;if(GCD(2*a-1-amm,2*b-1-bmm)>1||amm>9||bmm>9){printf(@
\pv@"Bad exponents\n");exit(1);}mw=m-2*(a+b-2);printf("\n");for(j=bmm-1;j>0;j--){for(i=1;i<amm;i++){if@
\pv@(Cn[i][j]>0){printf("%d,%d",i,j);s(2);}else s(5);}printf("\n");}for(j=bmm-1;j>0;j--){for(i=1;i<amm;@
\pv@i++){if(Cn[i][j]>0)prn(Cw[i][j]);else s(3);}printf("\n");}for(i=1;i<=nump;i++){for(j=1;j<=nump;j++)@
\pv@{if(J[i-1][1]<=J[j-1][1]&&J[i-1][2]<=J[j-1][2])O[i][j]=1;else O[i][j]=0;}}Nn[0]=1;for(i=1;i<=nump;@
\pv@i++)S[0][0][i]=1;for(k=1;k<nump;k++){n=0;for(e=0;e<Nn[k-1];e++){for(i=1;i<=nump;i++){for(j=1;j<=@
\pv@nump;j++)R[j]=S[k-1][e][j];if(R[i]!=0){R[i]=0;fl=0;for(p=1;p<nump&&fl==0;p++){if(R[p]!=0){for(q=@
\pv@p+1;q<=nump&&fl==0;q++){if(R[q]==0&&O[p][q]==1)fl=1;}}}if(fl==0){for(l=1;l<=mw;l++)SG[l]=0;SG[0]@
\pv@=1;for(q=1;q<=nump;q++){if(R[q]!=0){for(l=0;l<mw;l++){if(SG[l]!=0){for(c=1;l+c*J[q-1][0]<=mw;c++)@
\pv@SG[l+c*J[q-1][0]]=1;}}}}for(l=0;l<=mw;l++)WW[l]=0;for(p=1;p<=nump;p++){if(R[p]!=0)WW[J[p-1][0]]=1@
\pv@;}for(l=1;l<=mw;l++){if(SG[l]!=0&&W[l]!=0&&WW[l]==0)fl=1;}}if(fl==0){for(p=0;p<n&&fl==0;p++){for@
\pv@(l=1;R[l]==S[k][p][l]&&l<=nump;l++){;}if(l>nump)fl++;}if(fl==0){for(l=1;l<=nump;l++)S[k][n][l]=@
\pv@R[l];n++;}}}}}Nn[k]=n;}n=0;for(k=1;k<nump;k++){for(l=0;l<Nn[k];l++){for(j=bmm-1;j>0;j--){for(i=1;@
\pv@i<amm;i++){if(Cn[i][j]>0&&S[k][l][Cn[i][j]]!=0)prn(Cw[i][j]);else s(3);}printf("\n");}n++;}}printf@
\pv@("There are %d nonempty bistable subsets:\n 1",n+1);p=1;for(k=1;k<nump;k++){p=p+Nn[k];printf("+%d"@
\pv@,Nn[k]);}printf("=%d\n",p);}void prn(int r){if(r>=10)s(1);else s(2);printf("%d",r);}void s(int r)@
\pv@{int i;for(i=0;i<r;i++)printf(" ");}@}

\sk
{\smaller\smaller RIMS Kyoto University, Kyoto 606-8502 Japan}
\end{document}